\documentclass[letter,11pt]{amsart}
\usepackage{amssymb, mathrsfs, tikz} 
\usepackage{graphicx}
\usepackage{color}
\usepackage[font=footnotesize]{caption}
\usepackage{mathtools}
\usepackage{pinlabel}
\usepackage{graphicx}
\graphicspath{ {images/} }
\usepackage{subfig}

\usepackage{hyperref}
\usepackage{float}
\usepackage{calc}

\usetikzlibrary{decorations.markings, patterns}
\tikzset{->-/.style={decoration={
  markings,
  mark=at position #1 with {\arrow{>}}},postaction={decorate}}}
  
  \dedicatory{Dedicated to Dennis Sullivan on the occasion of
  his 80th birthday}  
  
\def\thetitle{{Random walks on mapping class groups}}
\hypersetup{
pdftitle=  \thetitle,
pdfauthor=  {Hyungryul Baik, Inhyeok Choi}
}

\newtheorem{thm}{Theorem}[section]
\newtheorem{lem}[thm]{Lemma}

\newtheorem{prop}[thm]{Proposition}

\newtheorem{defn}[thm]{Definition}

\theoremstyle{remark}

\newtheorem*{rem}{Remark}
\usepackage{tikz}
\usetikzlibrary{calc,decorations.markings}
\usetikzlibrary{shapes,snakes}

\theoremstyle{definition}
\newtheorem*{defn*}{Definition}
\newcommand{\nc}{\newcommand}
\nc{\dmo}{\DeclareMathOperator}
\dmo{\ra}{\rightarrow}
\dmo{\Prob}{\mathbb{P}}
\dmo{\E}{\mathbb{E}}
\dmo{\N}{\mathbb{N}}
\dmo{\Z}{\mathbb{Z}}
\dmo{\Q}{\mathbb{Q}}
\dmo{\R}{\mathbb{R}}
\dmo{\C}{\mathcal{C}}
\dmo{\X}{\mathcal{X}}
\dmo{\U}{\mathcal{U}}
\dmo{\T}{\mathcal{T}}
\dmo{\F}{\mathcal{F}}
\dmo{\AC}{\mathcal{AC}}
\dmo{\w}{\omega}
\dmo{\MIN}{\mathcal{MIN}}
\dmo{\Mod}{Mod}
\dmo{\PMod}{PMod}
\dmo{\PMF}{\mathcal{PMF}}
\dmo{\Mat}{Mat}
\dmo{\supp}{supp}
\dmo{\UE}{\mathcal{UE}}
\dmo{\vol}{vol}
\dmo{\B}{B}
\dmo{\PB}{PB}
\dmo{\PR}{PSL(2,\mathbb{R})}
\dmo{\GL}{GL(k, \mathbb{C})}
\dmo{\SL}{SL(2, \mathbb{Z})}
\dmo{\Isom}{Isom}
\dmo{\RP}{\mathbb{R} \mathrm{P}}
\dmo{\I}{\mathcal{I}}
\dmo{\el}{\ell_{\C}}
\dmo{\NN}{\mathcal{N}}
\dmo{\rk}{rank}
\dmo{\tr}{tr}
\dmo{\llangle}{\langle\langle}
\dmo{\rrangle}{\rangle\rangle}
\dmo{\Unif}{Unif}
\dmo{\Out}{Out}

\usetikzlibrary{decorations.markings, patterns}
\tikzset{->-/.style={decoration={
  markings,
  mark=at position #1 with {\arrow{>}}},postaction={decorate}}}

\numberwithin{equation}{section}

\title[random walks on mapping class groups]\thetitle

\author{Hyungryul Baik}
\address{Department of Mathematical Sciences, KAIST,  
291 Daehak-ro, Yuseong-gu, Daejeon 34141, South Korea }
\email{hrbaik@kaist.ac.kr}

\author{Inhyeok Choi}
\address{Department of Mathematical Sciences, KAIST,  
291 Daehak-ro, Yuseong-gu, Daejeon 34141, South Korea }
\email{inhyeokchoi@kaist.ac.kr}

\begin{document}
\maketitle

\begin{abstract}
This survey is concerned with random walks on mapping class groups. We illustrate how the actions of mapping class groups on Teichm{\"u}ller spaces or curve complexes reveal the nature of random walks, and vice versa. Our emphasis is on the analogues of classical theorems, including laws of large numbers and central limit theorems, and the properties of harmonic measures, under optimal moment conditions. We also explain the geometric analogy between Gromov hyperbolic spaces and Teichm{\"u}ller spaces that has been used to copy the properties of random walks from one to the other.
\end{abstract}



%
%

\section{Introduction} \label{sec:intro} 

A random walk is an example of a Markov chain or more generally a stochastic process.  Various models of random walks have been suggested and applied to settings including physics, economics, finance, biology, ecology, and many others. Among them are random walks on groups that deal with products of random group elements chosen by the same probability distribution. Since Kesten's pioneering work \cite{kesten1959symmetric}, the connection between the group structure (e.g. solvability, amenability, etc.) and the asymptotic behavior of random walks on the group has been studied in depth. This turned out to be fruitful in the perspective of both probability theory and group theory, and even geometry since the relatively new field called geometric group theory has arisen. 

In this survey, we study random walks on mapping class groups via their actions on Teichm{\"u}ller spaces. These actions are generalizations of the action of $\SL$ on $\mathbb{H}^{2}$. Here, $\SL$ can be viewed from two important perspectives in the study of non-commutative random walks: one as a discrete group acting on a negatively curved space, and the other one as a lattice in a Lie group that acts on a homogeneous space. Kaimanovich and Masur relate random walks on mapping class groups with these perspectives, which we explain in Section \ref{section:finiteness}.

The prototypes in the first perspective are random walks on free groups. In this case, random walks escape to infinity and reveal harmonic properties of the group. Explicit computations regarding this escape led to the analogues of classical limit laws including laws of large numbers, central limit theorems, and local limit theorems. For these computations, we refer the readers to \cite{woess2000random}, \cite{ledrappier2001free} and the references therein.

Appropriate notions that generalize free groups are hyperbolic groups introduced by M. Gromov. The corresponding geometric property is $\delta$-hyperbolicity: hyperbolic groups admit a geometric action on a proper $\delta$-hyperbolic space. Thanks to recent developments, one can deal with not only random walks on hyperbolic groups, but also those on weakly hyperbolic groups acting on non-proper $\delta$-hyperbolic spaces. We discuss Kaimanovich's theory for random walks on hyperbolic groups and Maher-Tiozzo's theory for random walks on weakly hyperbolic groups in Section \ref{section:curve}.

However, Teichm{\"u}ller spaces are not $\delta$-hyperbolic and mapping class groups are not lattices in semi-simple Lie groups of higher rank. Hence, both perspectives do not apply directly to mapping class groups (as Kaimanovich and Masur explain). Despite this contrast, we will first pursue the former perspective and make a due modification. Namely, Masur-Minsky's theory guided the usage of the curve complex, a non-proper $\delta$-hyperbolic space, in the study of the mapping class group. We explain how to compare the actions of mapping class groups on  Teichm{\"u}ller spaces and curve complexes in Section \ref{section:Teich}.

Having established relevant theories, Section \ref{section:displacement} and Section \ref{section:translation} deals with limit theorems for random walks on mapping class groups. We finish this survey by explaining the counting problem in mapping class groups and suggesting future directions.

Our survey is certainly not exhaustive at all. Especially, we regret not to explain the ideas of Sisto \cite{sisto2018contract}, Arzhantseva-Cashen-Tao \cite{arzhantseva2015growth} and Yang \cite{yang2019statistic}, \cite{yang2020genericity} that make use of contracting elements. This idea naturally covers the case of hyperbolic groups, CAT(0)-groups, mapping class groups and right-angled Artin groups. We also lack the explanation on Martin boundaries and its comparison with other boundaries.

Let us remark the connection between the theory of random walks and Patterson-Sullivan theory. Given a negatively curved manifold $X = \tilde{X} / \Gamma$ and its universal cover $\tilde{X}$, both theories construct measures on $\partial \tilde{X}$ using the deck transformation of $\tilde{X}$. Let us first consider a random walk on $\Gamma$ generated by the transition probability $\mu$. By applying this random walk to a point $x \in \tilde{X}$, one obtains a $\mu$-harmonic measure $\nu_{x}$ as the weak limit of the sample distribution at step $n$. In Patterson-Sullivan theory, each orbit point $gx$ is assigned the mass $e^{-s d(x, gx)}$, which sums up to a measure $\nu'_{x, s}$. By taking the weak limit of $\nu'_{x, s}$ as $s$ approaches to the growth exponent $\delta$, we obtain the Patterson-Sullivan measure $\nu'_{x}$ that is $\Gamma$-conformal with dimension $\delta$. These measures entail rich information about the geometry of $X$ (number of loops, etc.) and the dynamics on $X$ (mixing geodesic flow, etc.). Furthermore, in certain circumstances, $\tilde{X}$ is symmetric if and only if $\nu_{x}$ and $\nu'_{x}$ are proportional. 

Given these observations, the next goal is to build a parallel theory for Teichm{\"u}ller spaces. Harmonic measures and conformal densities will shed light on the geometry of the moduli space, and it matters to clarify which measure serves which role. In particular, given the non-homogeneity of Teichm{\"u}ller spaces, one can ask whether two measures on $\PMF$ differ or not.

Another important goal is to implement Patterson-Sullivan theory on groups, which dates back to Coornaert's work \cite{coornaert1993mesures}. In \cite{gekhtman2018counting}, Gekhtman-Taylor-Tiozzo utilized the automatic structure of the group to decompose the Patterson-Sullivan measure into countably many harmonic measures. This is intimately related to the counting problem on the group, yet the theory has not been extended to mapping class groups in full generality.

\subsection{Acknowledgements}
We thank Camille Horbez, Dongryul M. Kim, Hidetoshi Masai, \c{C}a\u{g}ri Sert for helpful discussions. The first author was partially supported by the National Research Foundation of Korea(NRF) grant funded by the Korea government(MSIT) \\ (No.2020R1C1C1A01006912) and the second author was supported by Samsung Science \& Technology Foundation grant No. SSTF-BA1702-01.

%
%
\section{Preliminary} \label{sec:prelim}

We review preliminary knowledge and fix conventions. Unless stated otherwise, $G$ denotes a finitely generated group that acts on a metric space $X$ by isometries. $S$ denotes a finite generating set of $G$. All measures are probability measures, and $\mu$ always denote a non-elementary measure on $G$. `Almost every' and `almost surely' are abbreviated to `a.e.' and `a.s.', respectively.

Given $\mu$, we can consider the step space $(G^{\Z}, \mu^{\Z})$, the product space of $G$ equipped with the product measure of $\mu$. The random walk on $G$ generated by $\mu$ is constructed by assigning to each step path $(g_{n})$ the sample path \[
\w_{n} := \left\{\begin{array}{cc} g_{1} \cdots g_{n} & n > 0 \\ id & n = 0 \\ g_{0}^{-1} \cdots g_{n+1}^{-1} & n < 0. \end{array}\right.
\]

\begin{rem}
It is worth making distinction between random walks growing from the right and from the left. Random walks in our convention grows from the right. Random walks on the ambient space $X$ are then modeled by applying these random walks to a reference point $o \in X$. One advantage of this model is that a.e. orbit path converges to a boundary point.

Meanwhile, random walks growing from the left naturally arise when we successively apply random isometries on an object (see \cite{horbez2018central} and \cite{benoist2016central}). Although the asymptotic behavior of orbit paths differs in two models, the limit theorems for displacements and translation lengths can be copied from one setting to the other one via the inversion \[
\w_{n} = g_{1} \cdots g_{n} \leftrightarrow \w_{n}' = g_{n}^{-1} \cdots g_{1}^{-1}
\]
and the identities $d(o, \w_{n}o) = d(o, w_{n}' o)$, $\tau(\w_{n}) = \tau(\w_{n}')$.
 \end{rem}
 
For $x, y, z \in X$, we denote the \emph{Gromov product} of $y, z$ with respect to $x$ by $(y, z)_{x}$. For details on the Gromov products, Gromov hyperbolic spaces and the Gromov boundary, see \cite{ghys1990gromov}, \cite{vaisala2005hyp}, \cite{bridson2013metric} and \cite{das2017gromov}.

We are interested in one group and two spaces associated to a closed orientable surface $\Sigma$ with genus at least 2: its mapping class group $\Mod(\Sigma) := \operatorname{Diff}^{+}(\Sigma) / \operatorname{Diff}_{0}(\Sigma)$, its Teichm{\"u}ller space $\T(\Sigma)$ and its curve complex $\C(\Sigma)$. The set of projective measured foliations on $\Sigma$ is denoted by $\PMF(\Sigma)$. $\MIN(\Sigma) \subseteq \PMF(\Sigma)$ denotes the set of projective measured foliations that correspond to minimal foliations, and $\widetilde{\MIN}(\Sigma)$ denotes its quotient by the equivalence relation of having trivial intersection. Finally, $\UE(\Sigma)$ denotes the set of uniquely ergodic foliations on $\Sigma$, which can be regarded as a subset of both $\MIN(\Sigma)$ and $\widetilde{\MIN}(\Sigma)$.

One dynamical quantity of an isometry $g$ of $X$ is the (asymptotic) translation length of $g$ defined by \[
\tau(g) := \lim_{n \rightarrow \infty} \frac{1}{n} d(o, g^{n} o).
\]
Recall that an isometry $g$ of a Gromov hyperbolic space $X$ falls into exactly one of the following categories: \begin{itemize}
\item $g$ has a bounded orbit (elliptic);
\item $g$ is not elliptic and has a unique fixed point in $\partial X$ (parabolic);
\item $g$ has two fixed points, an attractor and a repeller, in $\partial X$ (loxodromic).
\end{itemize}
In particular, $g$ is loxodromic if and only if $n \mapsto g^{n} x$ is a quasi-isometry for one (hence all) $x \in X$, if and only if $\tau(g) >0$.

As in the case of Gromov hyperbolic spaces, Thurston suggested $\PMF(\Sigma)$ as a natural boundary of $\T(\Sigma)$ in order to discuss the dynamical nature of mapping classes. Each $g \in \Mod(\Sigma)$ acts on $\T(\Sigma) \cup \PMF(\Sigma) \simeq \overline{B^{6g-6}}$ as a self-homeomorphism and has a fixed point. We have four cases: \begin{enumerate}
\item if $g$ has a fixed point in $\T(\Sigma)$, then $g$ is of finite order (periodic);
\item if $g$ fixes the projective class of a rational foliation, then $g$ fixes a multicurve (reducible);
\item if $g$ fixes an arational measured foliation without scaling, then $g$ is of finite order (periodic);
\item if $g$ scales up an arational measured foliation (hence fixes its projective class), then $g$ is pseudo-Anosov.
\end{enumerate}
Here, the final case and the other cases are mutually exclusive. In terms of their actions on $\C(\Sigma)$, periodic and reducible mapping classes are elliptic and pseudo-Anosov mapping classes are loxodromic \cite{masurMinsky2}.

%
%
\section{Early works and ergodic theorems} \label{section:finiteness}

In random walks, the step shift is ergodic and ergodic theory applies to the path space. Perhaps the first example of this strategy is that of Furstenberg and Kesten \cite{furstenberg1960random} for the product of random matrices. This is generalized to so-called Kingman's subadditive ergodic theorem. We present below one version of this theorem.

\begin{thm}[{\cite[Theorem 8.10]{woess2000random}}]\label{thm:Woess}
Let $(\Omega, \Prob)$ be a probability space and $U : \Omega \rightarrow \Omega$ be a measure-preserving transformation. If $W_{n}$ is a non-negative real-valued random variables on $\Omega$ satisfying the subadditivity $W_{n+m} \le W_{n} + W_{m} \circ U^{n}$ for all $m, n \in \N$, and $W_{1}$ has finite first moment, then there is a $U$-invariant random variable $W_{\infty}$ such that \[
\lim_{n \rightarrow \infty} \frac{1}{n} W_{n} = W_{\infty}
\]
almost surely and in $L^{1}(\Omega, \Prob)$. If $U$ is ergodic in addition, then $W_{\infty}$ is constant a.e.
\end{thm}

The subadditive ergodic theorem is particularly useful in non-commutative settings, where one cannot directly bring the results for Euclidean random walks. For example, let us consider a random walk $\w = (\w_{n})$ on $G$. After fixing a reference point $o \in X$, the displacement $d(o, \w_{n} o)$ becomes subadditive and the existence of the escape rate follows from the subadditive ergodic theorem. Additional properties of $G$ and its action on $X$, such as the non-amenability of $G$ and the properness of the action, guarantee that the escape rate is strictly positive and the random walk escapes to infinity. Hence, one can discuss the hitting measure induced on a suitable boundary of $X$, which is useful to investigate the asymptotics of the random walk. 

In fact, such boundary can be achieved \emph{a priori} without referring to the action of $G$ on $X$. Namely, Furstenberg extracted the ergodic component of the random walk with transition probability $\mu$ on a semi-simple Lie group $G$. This led to a suitable boundary $\partial G$ (called the \emph{Poisson boundary} of $(G, \mu$)) satisfying the Poisson property. It is however nontrivial to relate this measure-theoretical space with the space $X$ (in Furstenberg's case, $G$ acts on $X=G$ itself). Furstenberg utilized the structure of semi-simple Lie group to accomplish this, which was the first step for the proof of the following rigidity theorem.

\begin{thm}[{\cite[Theorem 1]{furstenberg1967random}}]\label{thm:furstenberg}
For $d \ge 2$ and $n \ge 3$, no countable group can become a cocompact lattice of $\Isom(\mathbb{H}^{d})$ and $\mathrm{SL}(n, \mathbb{R})$ simultaneously.
\end{thm}

Bringing this perspective to mapping class groups is due to Kaimanovich and Masur. After Thurston's compactification of $\T(\Sigma)$ with $\PMF(\Sigma)$ that led to the Nielsen-Thurston classification, Kaimanovich and Masur asked whether $\PMF(\Sigma)$ is the correct boundary of $\T(\Sigma)$ in the measure-theoretical viewpoint, that means, it hosts a $\mu$-stationary measure $\nu$ such that $(\PMF, \nu)$ becomes the Poisson boundary of $(\Mod, \mu)$. We present Kaimanovich-Masur's result below:

\begin{thm}[{\cite[Theorem 2.2.4, 2.3.1]{kaimanovich1996poisson}}]\label{thm:kaiMasur}
There exists a unique $\mu$-stationary $\nu$ on $\PMF(\Sigma)$, which is purely non-atomic and concentrated on $\UE \subseteq \PMF$, and $(\UE, \nu)$ is a $\mu$-boundary. In fact, $\nu$ is the hitting measure of $\mu$ on $\PMF$; for any $x \in \T(\Sigma)$ and $\Prob$-a.e. sample path $\w = (\w_{n})$, $\w_{n} x$ converges in $\PMF$ to a limit $F=F(\w) \in \UE$, and the distribution of the limits $F(\w)$ is given by $\nu$.

If $\mu$ has a finite entropy and finite first logarithmic moment with respect to the Teichm{\"u}ller metric in addition, then $(\PMF, \nu)$ is the Poisson boundary of $(\Mod, \mu)$.
\end{thm}

The proof relies on the geometry of $\T(\Sigma)$ and the structure of $\PMF$, which we briefly sketch now. As the escape to infinity is not established \emph{a priori}, the $\mu$-stationary measure $\nu$ is not constructed as the hitting measure of $\mu$; rather, $\nu$ is constructed indirectly and the escape to infinity follows from the property of $\nu$. The existence of $\nu$ relies on the fact that $\PMF(\Sigma)$ is a compact sphere. Using the structure of $\PMF$, namely, that $\PMF \setminus \MIN$ admits a countably infinite partition that respects $\Mod(\Sigma)$-action, one can deduce that $\nu$ is concentrated on $\MIN$. 

Let us now consider the quotient measure $\tilde{\nu}$ of $\nu$ on $\widetilde{\MIN}$, the set of equivalence classes of minimal foliations. Then a.e. path $(\w_{n})$ has the limit point $F(\w_{n}) \in \widetilde{\MIN}$ in the sense that $\w_{n} \tilde{\nu} \rightarrow \delta_{F(\w)}$ weakly. This map $\lambda : \w \mapsto \widetilde{\MIN}$ now serves as a Radon-Nikodym derivative for the $\mu$-boundary $(\widetilde{\MIN}, \tilde{\nu})$ and the Poisson boundary of $(\Mod, \mu)$ becomes nontrivial. This implies that the random walk escapes to infinity almost surely. The final technical step is to show that actually $\nu$ is concentrated on $\UE$ so that $\tilde{\nu}$ coincides with $\nu$. The uniqueness of $\nu$ now follows from the integral representation: $\lambda$ actually serves as a Radon-Nikodym derivative for $\nu$. Now the sample convergence $\w_{n} o \rightarrow F(\w)$ for a.e. $\w$ follows from the properties of uniquely ergodic foliations and universally convergent sequences.

In order to show that $(\PMF, \nu)$ is maximal, we invoke the strip approximation criterion introduced in \cite{kaimanovich2000hyp}. Namely, given that the entropy is finite, the maximality of $(\PMF, \nu)$ follows once we construct a measurable $\Mod(\Sigma)$-equivariant ``strips" $S : (F_{-}, F_{+}) \in \UE \times \UE \mapsto S(F_{-}, F_{+}) \subseteq \Mod(\Sigma)$ such that for all $g \in G$ and $\nu_{-} \otimes \nu_{+}$-a.e. $(F_{-}, F_{+}) \in \UE \times \UE$, \[
\frac{1}{n} \log \#(S(F_{-}, F_{+})g \cap B(e, |\w_{n}|)) \rightarrow 0 \quad \textrm{as}\quad n \rightarrow \infty
\]
in probability, where $|\w_{n}|$ denotes the word metric of $\w_{n}$ with respect to some finite generating set of $\Mod(\Sigma)$. Here we take $S(F_{-}, F_{+}) = \{h \in \Mod(\Sigma) : d(ho, [F_{-}, F_{+}]) \le M\}$ for some suitable $M$. As $\Mod(\Sigma)$ acts on $\T(\Sigma)$ properly discontinuously, balls of radius $M$ (with respect to the Teichm{\"u}ller metric) in $\T(\Sigma)$ contains finitely many translates $ho$ of $o$ and $\#(S(F_{-}, F_{+})g \cap B(e, k))$ grows at most linearly along $k$. Given this, the bottleneck of the growth of $\#(S(F_{-}, F_{+})g \cap B(e, |\w_{n}|))$ becomes the growth rate of $|\w_{n}|$ in probability. It turns out that finite logarithmic moment of $\mu$ suffices to control the overall growth in a subexponential manner.

Recall now Theorem \ref{thm:furstenberg} of Furstenberg: there, the dichotomy between two types of lattices was given as follows. If $G$ is a lattice in a semi-simple Lie group of rank at least 2, then there exist a measure $\mu$ with $\supp \mu = G$ and a number $\epsilon > 0$ such that the following holds. If $\mu$-harmonic functions $f_{1}$, $f_{2}$ on $G$ satisfy that (1) $0 \le f_{1}, f_{2} \le 1$ on $G$ and (2) $f_{1}(e), f_{2}(e) \ge 0.5 - \epsilon$, then (3) $\min[f_{1}(g), f_{2}(g)]$ does not tend to zero as $g \rightarrow \infty$. In contrast, if $G$ is a lattice in $\Isom(\mathbb{H}^{d})$, then for any $\mu$ and $\epsilon>0$ we can construct $\mu$-harmonic functions $f_{1}$, $f_{2}$ that satisfy (1), (2) but not (3). Kaimanovich and Masur similarly constructed such $\mu$-harmonic functions on $\Mod(\Sigma)$ and deduced that $\Mod(\Sigma)$ is not a lattice in semi-simple Lie groups.

This storyline already shows the interplay between the phenomena inside $\T(\Sigma)$ and the limiting phenomena at $\PMF(\Sigma)$. More specifically, (a) the escape to infinity and (b) the finite growth rate of strips are intimately related to (c) the characterization of the Poisson boundary. Here, a part of (c) helped establish (a), while (b) contributed to establishing another part of (c). Later we will see opposite situations, where (a) directly leads to a part of (c) and (c) helps establish a variant of (b).

A final remark is on the usage of Teichm{\"u}ller geometry. Examples that motivated Kaimanovich-Masur's work include the Gromov boundary of Gromov hyperbolic spaces that possess the visual measure and Patterson-Sullivan measure \cite{patterson1976limit}, \cite{sullivan1979density} of negatively curved manifolds. It is expected that the harmonic measure and these measures are mutually singular in compact manifolds with non-constant negative curvature \cite{ledrappier1990harmonic}, \cite{ledrappier1995dynamics}. Since the action of $\Mod(\Sigma)$ on $\T(\Sigma)$ has co-finite volume and  $\T(\Sigma)$ has variable curvature in some sense, we can also expect the mutual singularity of those measures. This topic will be revisited later.

The main obstacle to this problem was that both $\Mod(\Sigma)$ and $\T(\Sigma)$ are not Gromov hyperbolic. Nevertheless, Kaimanovich and Masur exploited partial hyperbolicity of the Teichm{\"u}ller space in order to establish the escape to infinity (see Subsection 1.5 of \cite{kaimanovich1996poisson} for example). This alludes to a unified storyline for random walks on mapping class groups and hyperbolic groups, which we will see in the subsequent sections.

%
%

\section{Curve complex} \label{section:curve}

Another space on which $\Mod(\Sigma)$ acts is the curve complex $\mathcal{C}(\Sigma)$. Introduced by Harvey \cite{harvey1981curve} as an analogy for the Tits building, this complex became a central object in the study of mapping class groups thanks to Masur-Minsky's theory (\cite{masurMinsky1}, \cite{masurMinsky2}). In particular, Masur and Minsky proved that: \begin{enumerate}
\item The shortest curve projection $\pi : \T(\Sigma) \rightarrow \mathcal{C}(\Sigma)$ sends geodesics to quasi-geodesics (with uniform constant),
\item $\mathcal{C}(\Sigma)$ is $\delta$-hyperbolic for some $\delta=\delta(\Sigma)$, and
\item $\T(\Sigma)$ and $\Mod(\Sigma)$ are weakly relatively hyperbolic in the sense of Farb.
\end{enumerate}
In fact, the constant $\delta(\Sigma)$ can be taken as $17$ for any surface $\Sigma$ (\cite{hensel2015slim}; see also \cite{aougab2013uniform}, \cite{bowditch2014uniform} and \cite{clay2014uniform}). We warn that mapping class groups are not relatively hyperbolic in general. Note also that $\T(\Sigma)$ and $\C(\Sigma)$ are not quasi-isometric, nor $\Mod(\Sigma)$ and $\C(\Sigma)$. Furthermore, the orbit map $[\phi] \mapsto \phi o$ from $\Mod(\Sigma)$ to $\T(\Sigma)$ is not a quasi-isometric embedding (but see \cite{eskin2017rank}).

As $\T(\Sigma)$ does, $\mathcal{C}(\Sigma)$ also captures the dynamics of each mapping class in $\Mod(\Sigma)$. Recall that pseudo-Anosov mapping classes have positive translation lengths on $\C(\Sigma)$ and $\T(\Sigma)$. Comparing these two translation lengths is an interesting question.

To investigate random walks on $\C(\Sigma)$, let us recall the work of Kaimanovich on the Poisson boundary of hyperbolic groups \cite{kaimanovich2000hyp}. Here, hyperbolic groups are finitely generated groups whose Cayley graph is Gromov hyperbolic with respect to a finite generating set. In this case, the given Cayley graph $\Gamma$ is locally finite and the Gromov boundary $\partial \Gamma$ is a compact metrizable space.

Given $\mu$ on a hyperbolic group $G$, one can obtain a $\mu$-stationary measure $\nu$ on $\partial \Gamma$ by the compactness argument. Kaimanovich also mentions the moment/entropy condition and the property on the action of $G$ on $\Gamma$ that the strip criterion requires. This property is enjoyed by hyperbolic groups and their boundaries, so $(\partial \Gamma, \nu)$ is indeed the Poisson boundary of $(\Gamma, \mu)$ under mild moment condition. Here the ambient space $\Gamma$ can be replaced with other \emph{locally compact} Gromov hyperbolic spaces $X$. Unfortunately, $\mathcal{C}(\Sigma)$ is not locally compact since every vertex have infinite valency. Being noncompact, it is hard to obtain certain objects as weak limits of measures on $\mathcal{C}(\Sigma) \cup \partial \mathcal{C}(\Sigma)$. 

Another criterion that identifies $\mu$-boundaries with the Poisson boundary is the ray approximation (also known as geodesic tracking) criterion. This asserts that if there exist nice rays $(\pi_{1}(\xi), \pi_{2}(\xi), \ldots)$ on $G$ for each $\xi \in B$ of a $\mu$-boundary $(B, \nu)$, then $(B, \nu)$ is the Poisson boundary of $(\Gamma, \mu)$. Here the niceness condition is that the ray projections are measurable and that for a.e. $\w = (\w_{n})$, if $\w$ heads to $\xi \in B$, then $\pi_{n}(\xi)$ and $\w_{n}$ deviate from each other sublinearly with respect to some gauge function $\mathcal{G}$.

When $\Gamma$ is properly discontinuously acting on a proper space $X$, then $|g|_{\mathcal{G}} := d(o, go)$ serves as a gauge. In this situation, if the progress $d_{X}(o, \w_{n} o)$ of the random walk is sublinear, then the trivial $\mu$-boundary becomes the Poisson boundary and $\Gamma$ becomes amenable. Conversely, random walks on non-amenable groups show positive escape rate. Nonetheless, this strategy is also not applicable to $\Mod(\Sigma)$ acting on $\mathcal{C}(\Sigma)$, as $\mathcal{C}(\Sigma)$ is not proper and each $\alpha \in \mathcal{C}(\Sigma)$ has infinite stabilizer. These limitations necessiated a radically different approach for $\mathcal{C}(\Sigma)$.

The breakthrough was made by Maher and Tiozzo. In \cite{maher2018random} they are concerned with weakly hyperbolic groups, the case where $X$ is a separable geodesic Gromov hyperbolic space and $G$ is a non-elementary subgroup of $\operatorname{Isom}(X)$. In this setting, they constructed the Poisson boundary for $(G, \mu)$ on the Gromov boundary of $X$. Before delving into their work, let us explore preceding observations by Maher.

Maher first deduced in \cite{maher2011random} the following: when $\supp \mu$ is a non-elementary subgroup of $\Mod(\Sigma)$, the probability that $\w_{n}$ is not pseudo-Anosov decreases to zero as $n \rightarrow \infty$. Overall, Maher mixes the probabilistic feature of harmonic measure and the group-theoretical structure of $\Mod(\Sigma)$ to achieve this. The first ingredient of Maher's argument is the relative conjugacy bounds of non-pseudo-Anosovs. Given this, it would suffice to show the transience of the set $R$ of elements whose relative conjugacy length is bounded by some constant (this would even lead to a stronger result in the almost sure sense). 

One possible approach is to use the fact that the harmonic measure $\nu(\partial R)$ of the limit set of $R$ dominates the probability of recurrence of $R$. Unfortunately, the limit set of $R$ is the entire $\PMF$ so this strategy fails. Instead, we consider a subset $R_{k}$ of $R$ that is contained in the union of centralizers of elements with word norm $1, \ldots, k$. If each of these centralizers has infinite copies in $\supp \mu$ ($\ast$), then a similar argument as in \cite{kaimanovich1996poisson} shows that the centralizers have harmonic measure 0. It is now argued that $\limsup_{n} \Prob(\w_{n} \in R \setminus R_{k})$ decreases down to zero as $k \rightarrow \infty$. The final ingredient is to ensure ($\ast$): it can be guaranteed by passing to a finite cover of $\Sigma$, if necessary.

Maher then proved in \cite{maher2010linear} that random walks show linear progress in the curve complex metric. Note that in contrast with $\Mod(\Sigma)$ with the word metric, $\mathcal{C}(\Sigma)$ is not locally compact and thus the standard non-amenability argument does not apply. Maher's idea was to construct a nested \emph{halfspaces (also known as shadows)} associated to each trajectory and show their abundance. This notion records partial useful history of the random walk, rather than merely recording the final product, in the form of stopping times. This strategy also helped improve the result of \cite{maher2011random} into the exponential decay of the probability of non-pseudo-Anosovs, given that the transition probability $\mu$ is finitely supported \cite{maher2012exp}.

Since these results were established based on the action on $\mathcal{C}(\Sigma)$, the key challenge was to remove the properness of $X$ and instead rely on the Gromov inequality among points only. Another important ingredient was that $\T(\Sigma)$ and $\mathcal{C}(\Sigma)$ partially shares the boundary structure, which at least coincide when we are concerned with $\mu$-boundaries. Recall again Theorem \ref{thm:kaiMasur}: $(\PMF, \nu)$ is a $\mu$-boundary concentrated on $\UE \subseteq \widetilde{\MIN}$. Here, $\UE$ is not only a subset of $\partial \T(\Sigma)$ but also a subset of $\partial \mathcal{C}(\Sigma) = \widetilde{\MIN}$. Thanks to this coincidence, one can use the $\mu$-boundary obtained from the dynamics on $\T(\Sigma)$ to investigate the asymptotic behavior of the random walk on $\C(\Sigma)$.

This effort culminated in Maher-Tiozzo's extensive work in \cite{maher2018random}. Recall that  $X$ is only assumed to be separable, geodesic and Gromov hyperbolic, and $G$ is a non-elementary countable subgroup of $\Isom(X)$. Not assuming that $X$ is proper, this setting includes $\Mod(\Sigma)$ acting on $\mathcal{C}(\Sigma)$ and $Out(F_{n})$ acting on the complex of free factors. Considering the previous works, we do not aim to obtain the $\mu$-boundary of $G$ directly from $X$ but from some other space, and extract the dynamic phenomenon from the coinciding boundary structure.

To this end, Maher and Tiozzo exploits the horofunction boundary $\bar{X}_{\infty}^{h}$, the usage of which is hinted by Calegari-Maher's work \cite{calegari2015scl}. Using the map $\phi : \bar{X}_{\infty}^{h} \rightarrow \partial X$ onto the Gromov boundary, one may\begin{enumerate}
\item push measures on $\bar{X}_{\infty}^{h}$ forward to $\partial X$, and
\item copy the (weak-$\ast$) convergence of measures on $\bar{X}_{\infty}^{h}$ to that on $\partial X$.
\end{enumerate}
As the horofunction compactification $\bar{X}^{h}$ of a separable metric space $X$ is always compact, the existence of a $\mu$-boundary $\nu$ on $\bar{X}^{h}$ and its concentration on the boundary $\bar{X}_{\infty}^{h}$ follow. Also deduced is a.e. convergence of measures $(\w_{n} \nu)$ on $\bar{X}_{\infty}^{h}$. All of these results can be pushed forward to $\partial X$; for example, the invariant measure $\nu$ on $\bar{X}^{h}$ gives rise to another measure $\tilde{\nu}$ on $\partial X$.

Meanwhile, given the convergence of $(\w_{n}\tilde{\nu})$ in $\partial X$, its convergence to an atom follows from the Gromov hyperbolicity. This leads to the fact that a.e. $\w = (\w_{n})$ has a limit point $\lambda = \lambda(\w) \in \partial X$ such that $(\w_{n} \tilde{\nu}) \rightarrow \delta_{\lambda}$ weakly. The final touch is to deduce the sample convergence from the weak convergence using shadows and the Gromov hyperbolicity; this is not available on $\bar{X}^{h}$ but on $X \cup \partial X$. It also follows that $\tilde{\nu}$ is non-atomic; if not, wandering of the maximum atom yields a contradiction with the assumption that $G$ is non-elementary.

Having established the relationship between the invariant, non-atomic measure $\tilde{\nu}$ with the escape to infinity of sample paths, a plenty of dynamical properties of the random walk can be obtained. These include positive escape rate, geodesic tracking and the linear growth of translation lengths. We summarize these below: 

\begin{thm}[cf. {\cite[Theorem 1.2, 1.3, 1.4]{maher2018random}}]\label{thm:maherTiozzo}
\begin{enumerate}
\item(Positive drift) There exists $L>0$ such that \[
\lim_{n\rightarrow \infty} \frac{d_{X}(o, \w_{n}o)}{n} = L \quad \textrm{a.s.}
\]
\item(Geodesic tracking) If $\mu$ has finite first moment, for a.e. $\w = (\w_{n})$, there exists a quasigeodesic ray $\gamma$ such that \[
\lim_{n \rightarrow \infty} \frac{d_{X}(\w_{n}o, \gamma)}{n} = 0.
\]
\item(Growth of translation lengths) There exists $L>0$ such that \[
\Prob(\tau(\w_{n}) \le Ln) \rightarrow 0 \,\,\textrm{as}\,\, n \rightarrow \infty.
\]
\end{enumerate}
\end{thm}

The final term is related to our main concern. Since a mapping class is pseudo-Anosov precisely when it acts on the curve complex loxodromically, we deduce that random walks on mapping class groups eventually become loxodromic in probability. 

A distinction should be made for random walks with bounded support. In this case, the argument of \cite{maher2012exp} indicates that the probability of shadows decrease exponentially as the distance from the origin increases. Using this, the above results are promoted into the following form.

\begin{thm}[cf. {\cite[Theorem 1.2, 1.3, 1.4]{maher2018random}}]\label{thm:maherTiozzoFinite}
Suppose that $\mu$ has bounded support. Then the following hold.
\begin{enumerate}
\item There exists $L, K>0$ and $0<c<1$ such that \[
\Prob(d_{X}(o, \w_{n} o) \le Ln) \le Kc^{n}
\]
for all $n$. 
\item For a.e. $\w = (\w_{n})$, there exists a quasigeodesic ray $\gamma$ such that \[
\limsup_{n \rightarrow \infty} \frac{d_{X}(\w_{n}o, \gamma)}{\log n} <\infty.
\]
\item There exists $L, K>0$ and $0<c<1$ such that \[
\Prob(\tau(\w_{n}) \le Ln) \le K c^{n}
\]
for all $n$.
\end{enumerate}
\end{thm}

In particular, the escape to infinity and the linear growth of translation lengths occur almost surely, rather than in probability; the geodesic tracking occur in a logarithmic manner, which is stronger than the sublinear one. In fact, combining Maher-Tiozzo's theory with Benoist-Quint's theory (to be explained later) yields the almost sure phenomena under finite second moment condition, as Dahmani and Horbez remark. We note that \cite{calegari2015scl}, \cite{sisto2018contract} and \cite{rivin2008walks} are also concerned with finitely supported random walks and deduce that the probability of non-pseudo-Anosov elements decay exponentially. 

We have yet to discuss the ultimate goal that $(\partial X, \tilde{\nu})$ is indeed the Poisson boundary of $(G, \mu)$. This requires a mild geometric condition on the action, namely, the acylindricity. The statement holds given that $G$ acts on $X$ acylindrically and $\mu$ has finite entropy and first moment. 

Maher-Tiozzo's work reveals that if abundant loxodromics are guaranteed, the coupling of the group structure and the space is not always required for investigating dynamical features of random walks. Note that Theorem \ref{thm:maherTiozzo} does not require the action of $G$ on $X$ to be cocompact or properly discontinuous (which would also restrict the geometry of $X$). Meanwhile, the Gromov hyperbolicity of $X$ plays a significant role throughout the argument. We also require $X$ to be separable in order to control the topology of the horofunction compactification. In the next section, we will examine how critical these conditions are.

\section{Teichm{\"u}ller spaces} \label{section:Teich}

We now discuss random walks on $\Mod(\Sigma)$ with respect to its action on $\T(\Sigma)$. Recall that the translation length $\tau_{\T(\Sigma)}(g)$ of a pseudo-Anosov mapping class $g$ with respect to the Teichm{\"u}ller metric and the stretch factor $\lambda(g)$ of $g$ have the relationship $\lambda(g) =  \log \tau_{\T(\Sigma)}(g)$. Hence, investigating the asymptotics of translation lengths on Teichm{\"u}ller spaces can reveal topological/dynamical properties of generic mapping classes. 

An immediate difficulty is that Teichm{\"u}ller spaces are not Gromov hyperbolic. To observe this, consider a geodesic triangle with vertices $o, T_{A}^{n} o, T_{B}^{-n} o$ for a point $o \in \T(\Sigma)$ and Dehn twists $T_{A}$, $T_{B}$ along disjoint curves $A$, $B$; the Hausdorff distance of $[T_{A}^{n} o, T_{B}^{-n}o]$ from $[o, T_{A}^{n} o] \cup [o, T_{B}^{-n} o]$ increases logarithmically \cite{masur1995teichmuller}. Another evidence is that the part of $\T(\Sigma)$ where a collection of disjoint curves $\{\gamma_{1}, \cdots, \gamma_{n}\}$ is pinched resembles a product space $\T(\Sigma \setminus \{\gamma_{1}, \cdots, \gamma_{n}\}) \times \prod_{i=1}^{n} \mathbb{H}^{2}$ \cite{minsky1996ext}.

Despite this failure, Teichm{\"u}ller spaces ($\Mod(\Sigma)$, resp.) share many aspects with negatively curved spaces (hyperbolic groups, resp.). For example, Margulis' work on the exponential growth of volumes and Deck transformation orbits of a negatively curved manifold has an analogy in the setting of Teichm{\"u}ller spaces and mapping class group orbits \cite{abem2012lattice}. The uniform exponential growth of hyperbolic groups is also copied onto mapping class groups by the work of \cite{anderson2007uniform}. In the same vein, many efforts have been put to copy the `thin triangle phenomenon' from Gromov hyperbolic geodesic spaces onto Teichm{\"u}ller spaces. 

Let us begin with Duchin's work on the geodesic tracking {\`a} la Kaimanovich. Kaimanovich suggested two criteria for determining the Poisson boundary of groups, the sublinear geodesic tracking and the strip approximation. Given that the strip approximation was effective enough to determine the Poisson boundary of mapping class groups, Kaimanovich asked whether the other criterion works, i.e., random walks on $\Mod(\Sigma)$ acting on $\T(\Sigma)$ show sublinear geodesic tracking. 

In fact, Kaimanovich-Masur's work already guides to the right candidate for the approximating geodesic. Namely, a.e. $\w =(\w_{n})$ possess the limit point $F(\w) \in \UE$ such that $\w_{n}o$ converges to $\UE$ in the sense of Thurston. Each geodesic $[o, \w_{n} o]$ is recorded at the initial point $o$ with the initial quadratic differential $\varphi_{n} \in QD_{o}$. Using Masur's comparison of Thurston and visual boundaries, it follows that $\varphi_{n} \rightarrow \varphi$ in $QD_{o}$ and the geodesic $\gamma=\gamma(\w)$ with the initial quadratic differential $\varphi$ converges to $F(\w)$.

As hinted before, the nuisance is the thin part of $\T(\Sigma)$. If $\gamma(\w)$ were always living inside a thick part of $\T(\Sigma)$, then one could apply the theory for Gromov hyperbolic spaces explained in Section \ref{section:curve}. Although $\gamma$ is approximated by geodesics connecting thick points, however, $\gamma$ may take a long excursion to the thin part of $\T(\Sigma)$. This led Duchin to focus on the phenomenon inside thick parts \cite{duchin2005teich}. More precisely, Duchin showed that when $\mu$ has finite first moment, a.e. $\w = (\w_{n})$ possesses a geodesic $\gamma: [0, \infty) \rightarrow \T(\Sigma)$ beginning from $o$ such that $\frac{1}{n}d(\w_{n}o, \gamma) 1_{K}( \gamma(d(o, \w_{n} o)))$ converges to zero. Here $K$ denotes a thick part of $\T(\Sigma)$.

Duchin's approach was to bring one particular property of ``thin triangles'' in Gromov hyperbolic spaces to some collection of triangles in $\T(\Sigma)$. In order for a random walk to be aligned along a geodesic, it is favored that consecutive orbits $\w_{n} o$ form a sort of ``highly obtuse triangles"; in such case, $d(\w_{n-k} o, \w_{n} o) + d(\w_{n} o, \w_{n+k} o)$ would be comparable to $d(\w_{n-k} o, \w_{n+k} o)$. Assuming such distance relations, we now conversely hope that each $\w_{n} o$ is not far away from the limiting geodesic $\gamma$. Motivated by this, Duchin required the following property. Let us first fix $A>0$, and consider a geodesic triangle $\triangle xyz$ with the longest side $[y, z]$. Let $w \in [y, z]$ be such that $d(x, y) = d(w, y)$. Then the desired property is \[
d(w, x) < A[d(y, x) + d(x, z) - d(y, z)].
\]
For example, $A=2$ works for triangles in $\mathbb{R}$-tree. In general, $A=2$ works for geodesic triangles in a Gromov hyperbolic space with a sidelength threshold. Duchin showed that geodesic triangles such that $w \in K$, together with a sidelength threshold, satisfy this property for some $A=A(K)$. The condition $w \in K$ led to the subsequence restriction in the theorem.

Before explaining how Rafi strengthened this approach, we digress to the complete sublinear geodesic tracking proved by Tiozzo \cite{tiozzo2015sublinear}. Tiozzo's approach is applicable not only to $\Mod(\Sigma)$ but also to groups acting on a proper Gromov hyperbolic spaces, groups with infinitely many ends, and groups acting on CAT(0) spaces. In the case of $\Mod(\Sigma)$ acting on $\T(\Sigma)$, we rely on the following fact: for a.e. $\w = (\w_{n})$, the forward limit $\eta$ and the backward limit $\xi$ are distinct points in $\UE$, hence transverse, and they are connected by a unique Teichm{\"u}ller geodesic. Given this, Tiozzo applies only the subadditive ergodic theorem to deduce the conclusion. 

Let us now discuss Rafi's analysis on thin triangles of $\T(\Sigma)$ in \cite{rafi2014hyperbolicity}. Motivated by the work of Masur-Minsky, Rafi aimed to investigate Teichm{\"u}ller geodesics with subsurface projection. Roughly speaking, a Teichm{\"u}ller geodesic $\gamma$ in $\T(\Sigma)$ for some surface $\Sigma$ can be cut into distinct subsegments $\gamma_{\alpha}$, each behaving like a Teichm{\"u}ller geodesic on some subsurface $Y_{\alpha}$ that is isolated during that time. Using this, Rafi deduced the following two instants of hyperbolicity in $\T(\Sigma)$.

The first item is fellow traveling. Consider two geodesics $\gamma : [a, b] \rightarrow X$ and $\eta : [a, b'] \rightarrow X$ with $d(\gamma(a), \eta(a)) < C$, $d(\gamma(b), \eta(b')) < C$. If $X$ is $\delta$-hyperbolic, then $\gamma$ and $\eta$ $K(C, \delta)$-fellow travel. We also expect $K(C, \epsilon)$-fellow traveling between such geodesics inside the $\epsilon$-thick part of $\T(\Sigma)$. However, there exists no \emph{a priori} uniform bound $K$ for every geodesics in $\T(\Sigma)$ having pairwise near endpoints. Rafi's theorem asserts that the geodesics $K(C, \epsilon)$-fellow travel if the pairwise near endpoints are $\epsilon$-thick, even if the geodesics are not entirely $\epsilon$-thick and visit the $\epsilon$-thin part.

The second item is as follows. Consider a geodesic triangle $\triangle xyz$ in $X$ and $p \in [yz]$. If $X$ is $\delta$-hyperbolic, then $p$ is within distance $K(\delta)$ from either $[x, y]$ or $[x, z]$. This is not guaranteed in $\T(\Sigma)$ in general, but there instead exist $K_{1}(\epsilon)$, $K_{2}(\epsilon)$ satisfying the following. If $p \subseteq I \subseteq [y, z]$ for some $\epsilon$-thick subsegment $I$ that is longer than $K_{1}(\epsilon)$, then the distance from $p$ and $[x, y] \cup [x, z]$ is at most $K_{2}(\epsilon)$. 

As we will see in the next section, these results are useful to compare the concatenation of geodesic segments $[x_{0}, x_{1}]$, $[x_{1}, x_{2}]$, $\ldots$, $[x_{N-1}, x_{N}]$ with the direct one $[x_{0}, x_{N}]$. However the intermediate journey during each segment is, if each segment behaves well near their endpoints, then the segments are aligned along $[x_{0}, x_{N}]$. This fact is exploited by Baik-Choi-Kim's pivoting that we explain later.

On the other hand, Rafi's approach that makes use of subsurface projections and marking distances was further exploited by Horbez, Dahmani-Horbez and Mathieu-Sisto.

Horbez's approach in \cite{horbez2018central} and Dahmani-Horbez's approach in \cite{dahmani2018spectral} begin with descending sample paths on $\T(\Sigma)$ to $\mathcal{C}(\Sigma)$ via the shortest curve projection $\pi : \T(\Sigma) \rightarrow \mathcal{C}(\Sigma)$ with some care. It is straightforward that the preimage of each point $p \in \mathcal{C}$ by $\pi$ is of infinite diameter. However, for a Teichm{\"u}ller geodesic $\gamma$ that is long in terms of both the Teichm{\"u}ller metric and the curve compelx metric, the (rough) preimage of $\pi(\gamma)$ may have stricter restriction. The following observation is motivated by the work of Dowdall-Duchin-Mausr improving Rafi's thin triangle result \cite[Theorem A]{dowdall2014statistical}.

\begin{prop}[{\cite[Proposition 3.7]{dahmani2018spectral}}]\label{prop:dahmaniHorbez}
For all $\kappa>0$, there exist $B, D>0$ such that the following holds. Let $[x, y]$ be a Teichm{\"u}ller geodesic that contains a subsegment $\gamma$ with sufficient progress on $\mathcal{C}(\Sigma)$, that means, $diam_{\mathcal{C}(\Sigma)}(\pi(\gamma)) > B$. If $z \in \T(\Sigma)$ satisfies that $\pi([x, z])$ crosses $\pi(\gamma)$ up to distance $\kappa$, then there exists a subsegment $\eta \subseteq [x, z]$ such that the Hausdorff distance of $\gamma$ and $\eta$ in $\T(\Sigma)$ is at most $D$ and $diam_{\mathcal{C}(\Sigma)}(\pi(\eta)) \ge diam_{\mathcal{C}(\Sigma)} (\pi(\gamma)) - B$. 
\end{prop}

In other words, the fellow-travelling among projections of long enough Teichm{\"u}ller geodesics can be lifted up. Recall also the result of Masur and Minsky that $\pi$ is coarsely $\Mod(\Sigma)$-equivariant, is coarsely Lipschitz, and sends Teichm{\"u}ller geodesics to $K(\Sigma)$-quasi-geodesics. Within this framework, we now explain how Dahmani and Horbez lifted the bahavior of random walks on $\C(\Sigma)$ to $\T(\Sigma)$.

Let us fix the reference point in $\T(\Sigma)$ by $o'$ temporarily, and let $o = \pi(o')$. We recall a result of Maher-Tiozzo: given a finitely supported measure $\mu$ on $\Mod(\Sigma)$, almost every sample path $\w = (\w_{n})$ of the random walk satisfies that \[
\lim_{n \rightarrow \infty} \frac{\tau_{\mathcal{C}(\Sigma)}(\w_{n})}{n} = \lambda,
\]
where $\lambda$ is the escape rate of the random walk in $\mathcal{C}(\Sigma)$. 

Let us consider geodesics $[o, \w_{n} o]$, $[\w_{n} o, \w_{n}^{2} o]$, $\ldots$, $[\w_{n}^{k-1} o, \w_{n}^{k}o]$ in $\mathcal{C}(\Sigma)$. Suppose, say, that $1000(\delta+B(10K'))\le d_{\mathcal{C}(\Sigma)}(o, \w_{n}o) \le 2\lambda n$ and \begin{equation}\label{eqn:weakDiscrep}
d_{\mathcal{C}(\Sigma)}(o, \w_{n} o) - \tau_{\mathcal{C}(\Sigma)}(\w_{n}) \le 0.01d_{\mathcal{C}(\Sigma)}(o, \w_{n} o)
\end{equation} (this will happen eventually in a.e. path $\w$). Then $[\w_{n}^{l}o, \w_{n}^{l-1} o]$ and $[\w_{n}^{l} o, \w_{n}^{l+1} o]$ should deviate early, at distance within $0.005 d_{\mathcal{C}(\Sigma)}(o, \w_{n} o)$. By $\delta$-hyperbolicity of $\mathcal{C}(\Sigma)$, there exist disjoint subsegments $[x_{l}, y_{l}]$ of $[o, \w_{n}^{k} o]$ that $K$-fellow travel with the middle $99\%$ of $[\w_{n}^{l-1} o, \w_{n}^{l} o]$. ($\ast$)

We now lift the situation to $\T(\Sigma)$ with the following ingredients. \begin{enumerate}
\item First, curve complex geodesics $[\w_{n}^{l-1} o, \w_{n}^{l} o]$ are close enough to the projections $\pi([\w_{n}^{l-1} o', \w_{n}^{l} o'])$ of the Teichm{\"u}ller geodesics, since the projections are quasi-geodesics and $\mathcal{C}(\Sigma)$ is $\delta$-hyperbolic. 
\item Similarly, $[o, \w_{n}^{k} o]$ and $\pi([o', \w_{n}^{k} o'])$ are close enough.
\item (1), (2) and ($\ast)$ imply that $\pi([o', \w_{n}^{k} o'])$ crosses the middle $98\%$ of each $\pi([\w_{n}^{l-1} o', \w_{n}^{l} o'])$ up to distance $K'$.
\end{enumerate}

We now apply Proposition \ref{prop:dahmaniHorbez} twice obtain subsegments $\eta_{l}$ of $[o', \w_{n}^{k} o]$ that satisfy the following. Let $\gamma_{0}$ be a subsegment of $[o', \w_{n} o']$ that projects onto the middle $96\%$ of $\pi([o', \w_{n}o'])$. Then $\eta_{l}$ and $\w_{n}^{l} \gamma_{0}$ are within Hausdorff distance $D(K)$ on $\T(\Sigma)$. Therefore, we have $d_{\T(\Sigma)}(\w_{n}^{l} o', [o', \w_{n}^{k} o']) \le d(o', \gamma_{0})$ for each $l$ and 
\[
d_{\T(\Sigma)}(o, \w_{n} o) - \tau_{\T(\Sigma)}(\w_{n}) = d_{\T(\Sigma)}(o, \w_{n} o) - \lim_{k} \frac{1}{k} d(o, \w_{n}^{k}o)\le 2 d_{\mathcal{T}(\Sigma)}(o', \gamma_{0}) + D(K)
\]
It now suffices to control the final term, the \emph{Teichm{\"u}ller} length of a left $2\%$ portion of $[o', \w_{n}o']$ with respect to the \emph{curve complex} distance. Although two distances are not comparable in general, the linear escape and sublinear tracking of a.e. sample path on both $\mathcal{C}(\Sigma)$ enables this. With this type of argument, Dahmani and Horbez obtains the following theorem: 

\begin{thm}[{\cite[Theorem 0.2]{dahmani2018spectral}}]\label{prop:dahmaniHorbezMain}
Suppose that $\mu$ is finitely supported. Then for a.e. sample path $(\w_{n})$, we have \[
\lim_{n \rightarrow \infty} \sqrt[n]{\lambda(\w_{n})} = \lambda,
\]
where $\log \lambda$ is the escape rate of the random walk.
\end{thm}

It is to be remarked that the finite support assumption originates from the spectral theorem for $\mathcal{C}(\Sigma)$. As Dahmani and Horbez explain, the arguments of Benoist-Quint and Maher-Tiozzo give rise to a spectral theorem for $\mathcal{C}(\Sigma)$ with finite second moment assumption. Given this, the rest of the Dahmani-Horbez's argument relies on the sublinear tracking and the subadditive ergodic theorem that only requires finite first moment.

Another way to relate the actions of $\Mod(\Sigma)$ on $\mathcal{C}(\Sigma)$ and $\T(\Sigma)$ was suggested by Mathieu and Sisto \cite{mathieu2020deviation}. Their philosophy is that nonelementary random walks on acylindrically hyperbolic groups are almost additive, so that most results can be reduced to that of commutative random walks on $\mathbb{R}$. For this purpose, they establish deviation inequalities, logarithmic geodesic tracking (see also \cite{sisto2017tracking}), and many more. These will be considered in the next section. 

Meanwhile, although it is true that $\Mod(\Sigma)$ is acting on $\mathcal{C}(\Sigma)$ and $\T(\Sigma)$ acylindrically, $\T(\Sigma)$ is not Gromov hyperbolic. Hence, one needs to bring the results on $\mathcal{C}(\Sigma)$ to $\T(\Sigma)$, which motivated Mathieu and Sisto to show the existence of $o \in \mathcal{C}(\Sigma)$ and $L \ge 0$ that satisfy the following. For $l_{1}, l_{2}, t \ge 0$ and $g, h \in \Mod(\Sigma)$ such that $d_{\C(\Sigma)}(go, ho) \ge L + l_{1} + l_{2}$, we have \[
diam^{\T(\Sigma)} \left[ \pi^{-1} (B^{\mathcal{C}(\Sigma)}(go, l_{1})) \cap N_{t}^{\T(\Sigma)} (\pi^{-1} (B^{\mathcal{C}(\Sigma)} (ho, l_{2})))\right] \le Lt,
\]
where $\pi$ denotes the shortest curve projection and $diam^{X}$, $N_{t}^{X}$, $B^{X}$ refer to the diameter, neighborhood and the ball with respect to $d_{X}$, respectively. This property follows from the coarse distance formula of the Teichm{\"u}ller metric in terms of (truncated) curve complex distances on subsurfaces \cite{rafi2007combi}, and bounded geodesic image theorem on curve complexes of subsurfaces with uniform constant (see \cite{masurMinsky1} and \cite{webb2013bounded}). Note that this property promotes the bounded distance of $\pi(p)$ from a long enough quasi-geodesic $\pi(\gamma)$ to the bounded distance of $p$ from $\gamma$.

Both Dahmani-Horbez's and Mathieu-Sisto's approach are concerned with geometric properties of $\pi$ that are not expected for arbitrary pairs of points or geodesics on $\T(\Sigma)$ but arise in almost every sample path. One partial reason, although not complete, is that $\Mod(\Sigma)$ acts on $\T(\Sigma)$ as isometries that translates $\epsilon$-thick reference point to another $\epsilon$-thick points, rather than arbitrary points. This implies that the randomness from random walks and other types of randomness in $\T(\Sigma)$ may show different behavior.

 In this spirit, Gadre-Maher-Tiozzo captured the contrast between the harmonic measure on $\PMF(\Sigma)$ arising from random walks and the Lebesgue measure \cite{gadre2017word}. A similar contrast holds between the Lebesgue measure on $\partial \mathbb{H}^{2}$ and the harmonic measure from a random walk on a cusped Fuchsian group. Gadre-Maher-Tiozzo considers the following quantity: for a boundary point $p \in \PMF(\Sigma)$, we first take a geodesic $\gamma$ tending to $p$, and approximate thick points $\gamma(t)$ with mapping class group orbits $h_{t} o$ of the reference point $o$. Then we compare the word norm of $h_{t}$ on $\Mod(\Sigma)$ and the displacement of $h_{t}$ with respect to the curve complex metric. In terms of the Lebesgue measure, $d_{\Mod(\Sigma)}(1, h_{t})$ grows indefinitely faster than $d_{\mathcal{C}(\Sigma)}(o, h_{t}o)$ in almost every choice of $p$; on the other hand, in terms of the harmonic measure for $\mu$ with finite first moment in the word metric on $\Mod(\Sigma)$, $d_{\Mod(\Sigma)}(1, h_{t})$ and $d_{\mathcal{C}(\Sigma)}(o, h_{t}o)$ are comparable and their ratio converges to a uniform constant in almost every choice of $p$. 

Let us finish this section by explaining a consequence of Rafi's theorems that will be used later on. Consider a geodesic triangle $\triangle xyz$ in $\T(\Sigma)$. A priori, $\triangle x y z$ is not $\delta$-thin and $[yz]$ need not be contained in a bounded neighborhood of $[xy] \cup [xz]$. However, suppose that $[x, y]$ initially fellow travels with a thick segment $[x, y']$. This forces that $[x, y]$ is initially thick also, and Rafi's theorem asserts that this beginning portion should be contained in a bounded neighborhood of $[x, y] \cup [y, z]$.

Let us similarly suppose that $[x, z]$ initially fellow travels with a thick segment $[x, z']$, and $[x, y']$ and $[x, z']$ are heading to different directions, i.e., $(y', z')_{x}$ is bounded. Then the initial segment of $[x, y]$ cannot be contained in the neighborhood of $[x, z]$ and vice versa. Finally, if we further suppose that points $y, z$ are also thick, then Rafi's fellow traveling theorem implies that $\triangle xyz$ is an obtuse thin triangle: $[y, z]$ and $[x, y] \cup [x, z]$ are within bounded Hausdorff distance.

\section{Limit theorems I: displacements} \label{section:displacement}

In the Euclidean setting, stronger moment assumptions lead to finer description on random walks. For example, it is believed that strong laws of large numbers (SLLN) are linked with the finitude of first moment; central limit theorems (CLT) and laws of the iterated logarithm (LIL) are relevant to the finitude of second moment; when the random walk has finite exponential moment, large deviation principles (LDP) is also available. Many recent work in this topic tried to bring these results to hyperbolic settings under suitable moment conditions. Among them we explain the results of Benoist-Quint, Horbez, Mathieu-Sisto, Boulanger-Matieu-Sert-Sisto, Gou{\"e}zel, Baik-Choi-Kim and Choi.

In hyperbolic settings, two meaningful quantities arise from random walks $\w = (\w_{n})$ on the isometry group: the displacement $d(o, \w_{n} o)$ of a reference point $o \in X$ and the translation length $\tau(\w_{n})$. The first one is subadditive while the later one is not; this complicates the investigation of translation lengths. We will first discuss the theorems for displacements and then move on to case of translation lengths.

The theorem in hyperbolic settings that corresponds to laws of large number is the subadditive ergodic theorem. For completeness, we spell out the statement: \begin{thm}
Let $X = \T(\Sigma)$ or $\C(\Sigma)$ and suppose that $\mu$ has finite first moment. Then there exists $\lambda >0$, called the escape rate of $\mu$, such that the random variables $\frac{1}{n} d_{X}(o, \w_{n} o)$ converge to $\lambda$ in $L^{1}$ and almost surely.
\end{thm}
Here the non-zero escape rate is due to the non-amenability of $\Mod(\Sigma)$ in the case of $X = \T(\Sigma)$, whereas it follows from the existence of `persistent joint' from Maher-Tiozzo's argument in the case of $X = \C(\Sigma)$ (which ultimately relies on the fact that the harmonic measure for $\mu$ on $\partial \C(\Sigma)$ is atom-free).

We remark that this is not a consequence of the Borel-Cantelli argument. Indeed, for example, the exponential decay of $\Prob(\frac{1}{n} d(o, \w_{n} o) \ge \lambda + \epsilon)$ for $\epsilon>0$ implies that $\mu$ has finite exponential moment. In contrast, $\Prob(\frac{1}{n} d(o, \w_{n} o) \le \lambda - \epsilon)$ does decay exponentially even without any moment condition due to the recent work of \cite{gouezel2021exp}. We will postpone the details of Gou{\"e}zel's technique at the moment; we note that the technique is powerful enough to deduce other results including the continuity of the escape rate.

The next natural goal is CLTs of the following form: \begin{thm}
Let $X = \T(\Sigma)$ or $\C(\Sigma)$. Suppose that $\mu$ is non-arithmetic and satisfies some moment condition. Then there exists $\sigma > 0$ such that $\frac{1}{\sqrt{n}}[d_{X}(o, \w_{n} o) - \lambda n]$ converges to the Gaussian law $\mathcal{N}(0, \sigma)$ in law, where $\lambda>0$ is the escape rate of the random walk. Explicitly, for any $a< b$, we have \[
\lim_{n \rightarrow \infty} \Prob\left[a \sqrt{n} \le d_{X}(o, \w_{n} o) - \lambda n \le b \sqrt{n}\right] = \int_{a}^{b} \frac{1}{\sqrt{2\pi} \sigma} e^{-x^{2}/2\sigma^{2}} \, dx.
\]
\end{thm}
This direction dates back to Sawyer-Steger's investigation \cite{sawyer1987free} on the random walks on free groups, which was also discussed by Ledrappier \cite{ledrappier2001free}. Its generalization to Gromov hyperbolic group under the finite exponential moment assumption is attributed to Bj{\"o}rklund \cite{bjorklund2010clt}. The current CLT under the finite second moment assumption was proven by Benoist and Quint in \cite{benoist2016central}, using the machinery of their previous work on linear groups. Finally, by the lifting principle that we explained before, Horbez generalized this CLT to Teichm{\"u}ller spaces \cite{horbez2018central}.

Benoist-Quint's setting is a non-elementary group $G$ acting on a proper, quasiconvex, Gromov hyperbolic space $X$ and a non-elementary, non-arithmetic Borel measure on $G$. Note that $G$ and $\mu$ need not be discrete here. The properness of $X$ is assumed to exploit the Gromov compactification $X \cup \partial X$ and the Busemann compactification $X \cup \partial_{B} X$. The main strategy is to find an alternative for the random variable $d(o, \w_{n} o)$, which can be expressed as martingales with step differences controlled in $L^{2}$ and in probability. The first trick is to use Busemann functions $\sigma(g, x) := \lim_{n \rightarrow \infty} [d(g^{-1} o, x_{n}) - d(o, x_{n})]$ for $x \in \partial_{B}X$ and $x_{n} \rightarrow x$ instead of displacements. In contrast with displacements, Buseman functions satisfy the cocycle condition\begin{equation}\label{eqn:additiveBusemann}
\sigma(gg', x) = \sigma(g, g'x) + \sigma(g', x).
\end{equation} At the cost of this advantage, however, one should pay attention to several details. First, the choice of the  boundary point $x$ causes some asymmetry. Moreover, the decomposition of $\sigma(\w_{n}, x)$ into $\sigma(g_{k+1}, \w_{k} x)$ as in Equation \ref{eqn:additiveBusemann} involves different boundary points $\w_{k}x$; hence, the argument requires analysis on the boundary action of $G$ and the stationary measure on $\partial_{B} X$ or $\partial X$. Second, $\sigma(\w_{n}, x)$ is nonetheless different from $d(o, \w_{n} o)$ and the discrepancy $d(o, \w_{n} o) - \sigma(\w_{n}, x)$ should be controlled for large enough $n$ in probability. Finally, the quantities $\sigma(g_{k+1}, \w_{k} x)$ still may not be adequate for martingale CLTs, as they are not `centered' at the right value, namely, the escape rate $\lambda$. One should therefore solve the cohomological equation and center $\sigma(g_{k+1}, \w_{k} x)$ by subtracting bounded random variables. After all these preliminary steps, one controls the step differences in $L^{2}$ and in probability using the finitude of second moment of $\mu$ and concludes the proof. Note again that the spirit of this proof is ``from the infinity", rather than ``working inside the space".

As we have seen before, the lifting argument often promotes phenomena in $\mathcal{C}(\Sigma)$ to the corresponding ones in $\mathcal{T}(\Sigma)$. Horbez's strategy in \cite{horbez2018central} was to lift the ingredients for Benoist-Quint's CLT, including the centrability of Busemann functions and the summable decay of shadows in particular directions, from $\mathcal{C}(\Sigma)$ to $\mathcal{T}(\Sigma)$. Although $\mathcal{C}(\Sigma)$ is not proper and thus Benoist-Quint's proof does not apply as is, theses ingredients are available on non-proper spaces by Maher-Tiozzo's work. Once the ingredients are lifted using Proposition \ref{prop:dahmaniHorbez}, Benoist-Quint's argument applies to Busemann functions on $\T(\Sigma)$ and the desired CLT follows.

Another approach to the CLT for displacements was proposed by Mathieu and Sisto in \cite{mathieu2020deviation}. In fact, they provide much more general framework, requiring the control on the defects of the form $Q_{n+m}(\w) - Q_{n}(\w) - Q_{m}(\theta_{n}\w)$ (the Gromov products $(\w_{n+m} o, o)_{\w_{n} o}$ in our setting, for example) and yielding quantitative estimations to which extent an addition with defects differ from the ideal addition.

To see the principle behind this, let $G$ be acting on any metric space $X$ and let us assume that $\E [(\w_{n} o, \check{\w}_{n} o)_{o}^{2}]$ is bounded by some constant $B$ for all $n$. We now estimate the distances among $o$, $\w_{n} o$, $\w_{2n} o$, $\ldots$, $\w_{2^{k} n} o$. If the points were always perfectly aligned, then $\E[ d(o, \w_{2^{k}n} o)]$ and $Var [d(o, \w_{2^{k}n} o)]$ would grow linearly with respect to $2^{k}$ (note that the family $\{\w_{in}^{-1} \w_{(i+1)n}\}_{i}$ consists of independent RVs). We would also have $d(o, \w_{2^{k} n} o) = \sum_{i=1}^{2^{k}} d(\w_{i-1} o, \w_{i} o)$ and the classical CLT would imply that $\frac{1}{\sqrt{2^{k}n}}[d(o, \w_{2^{k} n} o) - \E[d(o, \w_{2^{k} n} o)]]$ converges in law to a Gaussian law $\mathscr{N}(0, \sigma_{n})$, where $\sigma_{n} = \sqrt{Var[d(o, \w_{n} o)] / n}$. 

However, the addition does not happen perfectly in reality and the deficit is recorded in the form $2(\w_{in} o, \check{\w}_{kn} o)_{\w_{jn}o}$. Note that in order to make $d(o, \w_{2^{k} n} o)$ out of $d(o, \w_{n}o)$, $\ldots$, $d(\w_{(2^{k} - 1)n} o, \w_{2^{k} n} o)$, we need to add up $2^{k} - 1$ deficits $2(\w_{2^{t}(2j-2)} o, \w_{2^{t}(2j-1)} o)_{\w_{2^{t+1}j} o}$ for $t=0, \ldots, k-1$ and $j=1, \ldots, 2^{k-t-1} - 1$. Note that for a fixed $t$, these become a family of $(2^{k-t-1} - 1)$ i.i.d. with variance less than $B$. Summing them up and dividing by $1/\sqrt{2^{k} n}$, the error from these terms is bounded by $7/\sqrt[3]{n}$ outside an event of probability at most $8B/\sqrt[3]{n}$. By taking dyadic $n = 2^{m}$, we deduce that $\frac{1}{\sqrt{2^{m}}} [d(o, \w_{2^{m}} o) - \E[d(o, \w_{2^{m}}o)]]$ is Cauchy and $\sigma_{m} \rightarrow \sigma >0$ (here is required at least linear growth of $Var[d(o, \w_{2^{m}} o)]$, which is deduced from the non-arithmeticity of $\mu$). Similar argument can handle non-dyadic steps also, if $\E[(\w_{n} o, \check{\w}_{m} o)_{0}^{2}]$ is uniformly controlled for arbitrary $n$, $m$.

It remains to control $\E[(\w_{n} o, \check{\w}_{n} o)_{o}^{2}]$ as promised, for which Mathieu-Sisto's argument requires two assumptions: (1) that $\mu$ has finite exponential moment and (2) that the action is acylindrical. Although the second assumption can be removed due to the results in \cite{boulanger2020large}, it is enough for $\Mod(\Sigma)$ acting on $\C(\Sigma)$. Moreover, $\T(\Sigma)$ also fits into this scheme since it is acylindrically intermediate for $(\Mod(\Sigma), \C(\Sigma))$.

Let us explain how Mathieu-Sisto's viewpoint of `almost exact addition' was expanded in the works of Boulanger-Mathieu-Sert-Sisto, Gou{\"e}zel and Choi. All these works use the same modification of random walks as follows. Given a measure $\mu$ on $G$ and $S \subseteq G$ such that $\alpha := \min \{\mu(g) : g \in S\} > 0$, there exists a measure $\eta$ such that $\mu = \alpha \mu_{S} + (1-\alpha) \eta$, where $\mu_{S}$ is the uniform measure on $S$. We then consider: \begin{itemize}
\item Bernoulli RVs $\rho_{i}$ (with $\Prob(\rho_{i} = 1) = \alpha$ and $\Prob(\rho_{i} = 0) = 1-\alpha$), 
\item $\eta_{i}$ with the law $\eta$, and 
\item $\nu_{i}$ with the law $\mu_{S}$,
\end{itemize} 
all independent, and define\[
\begin{aligned}
\gamma_{i} = \left\{\begin{array}{cc}\eta_{i} & \textrm{when}\,\,\rho_{i} = 1, \\
\nu_{i} & \textrm{when}\,\,\rho_{i} = 0.\end{array}\right.
\end{aligned}
\]
Then $\gamma_{i}$ are i.i.d. of the law $\mu$, which models the random walk on $G$ generated by $\mu$. Let us also define $\mathcal{N}(k) := \sum_{i=1}^{k} \rho_{i}$ and $\vartheta(i) := \min\{j \ge 0 : B(j) = i\}$ for convenience. In this perspective, a random trajectory consists of relatively usual steps $(\gamma_{\vartheta(i) + 1}, \cdots, \gamma_{\vartheta(i+1) - 1})$ and special steps $\gamma_{\vartheta(i)}$ in an alternating way, the first one being chosen with law $\nu$ and the second one being chosen with law $\mu_{S}$. 

Morally, $\nu$ is designed to behave almost like $\mu$: they share the same moment condition and similar moment values. The displacements made by these usual steps are then linked with the special steps from $S$. The desired property of special steps is that they `almost align' consecutive displacements with high probability. This is encoded in the notion of \emph{Schottky set}, which stems from the classical Schottky decomposition. We note that the exact definitions of Schottky sets differ in the aforementioned three references; among them we introduce a version that works in all of three settings.

\begin{defn}\label{dfn:Schottky}
Let $K, K', \epsilon> 0$. A finite set $S$ of isometries of $X$ is said to be \emph{$(K, K')$-Schottky} if the following hold: \begin{enumerate}
\item for all $x, y \in X$, $|\{s \in S : (x, s^{i}y)_{o} \ge K\,\,\textrm{for some}\,\,i>0\}| \le2$;
\item for all $x, y \in X$, $|\{s \in S : (x, s^{i}y)_{o} \ge K\,\,\textrm{for some}\,\,i<0\}| \le2$;
\item for all $s \in S$ and $i \neq 0$, $d(o, s^{i}o) \ge K'$.
\end{enumerate}
When $X$ is Teichm{\"u}ller space, $S$ is said to be \emph{$(K, K', \epsilon)$-Schottky} if the following condition holds in addition to the above three: \begin{enumerate}\setcounter{enumi}{3}
\item for all $s \in S$ and $i \in \Z$, the geodesic $[o, s^{i}o]$ is $\epsilon$-thick.
\end{enumerate}
\end{defn}

By employing Schottky sets for `linking steps', one can add up step distances almost exactly. In particular, Boulanger-Mathieu-Sert-Sisto recovered the following deviation inequality of Mathieu-Sisto without the acylindricality assumption: if $\mu$ has finite exponential moment, there exists $K, \kappa >0$ such that \begin{equation}\label{eqn:mathieuSisto}
\Prob[(o, \w_{n}o)_{\w_{i} o} \le R] \le Ke^{-\kappa R}
\end{equation}
holds for any $0 \le i \le n$. This result is subsequently used to establish the following large deviation principle for the random walk:
\begin{thm}[{\cite[Theorem 1.1]{boulanger2020large}}]
If $\mu$ has finite exponential moment, then there exists a proper convex function $I : \R \rightarrow [0, \infty]$ (called the \emph{rate function}) that satisfies\[
-\inf_{\alpha \in \operatorname{int}(R)} \le \liminf_{n} \frac{1}{n} \ln \Prob\left[ \frac{1}{n} d(o, \w_{n} o) \in R\right] \le  \limsup_{n} \frac{1}{n} \ln \Prob\left[ \frac{1}{n} d(o, \w_{n} o) \in R\right] \le -\inf_{\alpha \in \bar{R}} I(\alpha)
\]
for any measurable subset $R$ of $\mathbb{R}$. Moreover, $I$ vanishes only at the escape rate $\lambda$ of the random walk.
\end{thm}
Roughly speaking, the probability that $\frac{1}{n} d(o, \w_{n}o)$ deviates from $\lambda$ decays exponentially, the speed of which is precisely encoded in $I$. We note that Boulanger-Mathieu-Sert-Sisto establishes the rate function from above (for values greater than $\lambda$) in much more general setting, on arbitrary metric spaces. Indeed, the existence of the rate function from above is essentially equivalent to the finitude of exponential moment, rather than the geometric property of the underlying space, as mentioned before. Meanwhile, establishing the rate function from below requires the existence of Schottky sets and the Gromov inequalities among points.

It was unexpected, however, that the exponential decay of the deviation from below does not require any moment condition.
\begin{thm}[{\cite[Theorem 1.1, 1.2]{gouezel2021exp}}]
\begin{enumerate}
\item If $\mu$ has finite first moment and $\lambda$ is its escape rate, then $\Prob[d(o, \w_{n} o) \le (\lambda - \epsilon)n]$ decays exponentially for any $\epsilon > 0$.
\item If $\mu$ has infinite first moment, then there is no finite `escape rate': $\Prob[d(o, \w_{n} o) \le rn]$ decays exponentially for any $r > 0$.
\end{enumerate}
\end{thm}

To establish this result, Gou{\"e}zel first takes suitable integer $N$ and a Schottky set $S$ so that we have a decomposition \begin{equation}\label{eqn:schottkyMeas}
\mu^{\ast N} = \alpha \mu_{S}^{\ast 2} + (1-\alpha)\nu.
\end{equation} Here the $N$-th convolution of $\mu$ is designed to guarantee sufficiently large size of $S$; the purpose of the self-convolution of the Schottky measure will become apparent soon. Then the composition $\gamma_{i}$ of the Bernoulli variable $\rho_{i}$ and $\eta_{i}$, $\nu_{i}$ models the convolution of steps $g_{N(i-1)+1}, \ldots, g_{Ni}$. We now record some of the special steps $\vartheta(i)$ as \emph{pivotal times}, which are meant to be the crucial moments throughout the history of the random path.

At each pivotal step, we hope that two Schottky segments are directed away from each other and the former (latter, resp.) Schottky segment does not cancel out the previous (upcoming, resp.) progress. To be concrete, suppose that we have chosen $\{m_{1}<\ldots<m_{k}\}$ from the special steps $\{\vartheta(i)\}_{i=1}^{M-1}$ as pivotal times for the path $(g_{1}, \ldots, g_{\vartheta(M)})$. Let
\[\begin{aligned}
w_{i} := \w_{N(m_{i} + 1)}^{-1} \w_{Nm_{i+1}}, \,\,s_{i} := \w_{Nm_{i}}^{-1} \w_{N(m_{i} +0.5)}, \,\,s_{i}' := \w_{N(m_{i}+0.5)}^{-1} \w_{N(m_{i} +1)}.
\end{aligned}
\]
and $w_{0} := \w_{Nm_{1}}$, $s_{0} = id$. The desirable situation that $m_{1}, \ldots, m_{k}$ already satisfy is the following for suitable $K$: \begin{itemize}
\item $(s_{i}^{-1} o, s'_{i}o)_{o} \le K$ for $i=1, \ldots, k$,
\item $(s'^{-1}_{i}o, w_{i}o)_{o} \le K$ for $i=1, \ldots, k$,
\item $(w_{i-1}^{-1} s'^{-1}_{i-1}o, s_{i} o)_{o} \le K$ $i=1, \ldots, k$. 
\end{itemize}

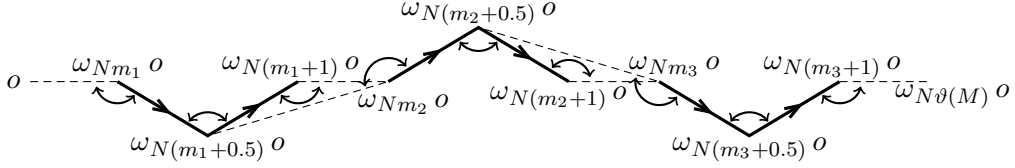
\begin{figure}[H]
\begin{tikzpicture}
\def\a{0.8}
\def\b{3}
\def\c{1.2}

\draw[very thick, decoration={markings, mark=at position 0.28 with {\draw (-0.2, 0.07) -- (0, 0) -- (-0.2, -0.07);}, mark=at position 0.78 with {\draw (-0.2, 0.07) -- (0, 0) -- (-0.2, -0.07);}}, postaction={decorate}]  (-\c, \c*0.6) -- (0, 0) -- (\c, \c*0.6);

\draw[very thick, shift={(\c*3, 0)}, decoration={markings, mark=at position 0.28 with {\draw (-0.2, 0.07) -- (0, 0) -- (-0.2, -0.07);}, mark=at position 0.78 with {\draw (-0.2, 0.07) -- (0, 0) -- (-0.2, -0.07);}}, postaction={decorate}]  (-\c, \c*0.6) -- (0, \c*1.2) -- (\c, \c*0.6);

\draw[very thick, shift={(\c*6, 0)}, decoration={markings, mark=at position 0.28 with {\draw (-0.2, 0.07) -- (0, 0) -- (-0.2, -0.07);}, mark=at position 0.78 with {\draw (-0.2, 0.07) -- (0, 0) -- (-0.2, -0.07);}}, postaction={decorate}]  (-\c, \c*0.6) -- (0, 0) -- (\c, \c*0.6);

\draw[densely dashed]  (\c, \c*0.6) --  (\c*2, \c*0.6) -- (0,0);
\draw[densely dashed]  (\c*4, \c*0.6) --  (\c*5, \c*0.6) -- (\c*3, \c*1.2);
\draw[densely dashed] (-\c, \c*0.6) -- (-2*\c, \c*0.6);
\draw[densely dashed] (7*\c, \c*0.6) -- (8*\c, \c*0.6);

\draw[thick, <->, shift={(-\c, \c*0.6)}, rotate=193] (0.24*\c, 0) arc (0:129:0.24*\c);
\draw[thick, <->, shift={(\c, \c*0.6)}, rotate=-13] (0.24*\c, 0) arc (0:-129:0.24*\c);
\draw[thick, <->, shift={(2*\c, \c*0.6)}, rotate=40] (0.26*\c, 0) arc (0:152:0.26*\c);
\draw[thick, <->, shift={(4*\c, \c*0.6)}, rotate=12] (0.24*\c, 0) arc (0:129:0.24*\c);
\draw[thick, <->, shift={(5*\c, \c*0.6)}, rotate=170] (0.26*\c, 0) arc (0:150:0.26*\c);
\draw[thick, <->, shift={(7*\c, \c*0.6)}, rotate=-13] (0.24*\c, 0) arc (0:-129:0.24*\c);

\draw[thick, <->, shift={(0,0)}, rotate=37] (0.26*\c, 0) arc (0:106:0.26*\c);
\draw[thick, <->, shift={(3*\c, 1.2*\c)}, rotate=217] (0.26*\c, 0) arc (0:106:0.26*\c);
\draw[thick, <->, shift={(6*\c,0)}, rotate=37] (0.26*\c, 0) arc (0:106:0.26*\c);

\draw (-2*\c - 0.18, \c*0.6) node {$o$};
\draw (-\c, \c*0.6 + 0.18) node {$\w_{Nm_{1}} o$};
\draw (0, -0.18) node {$\w_{N(m_{1} + 0.5)}o$};
\draw (\c-0.1, \c*0.6 + 0.18) node {$\w_{N(m_{1} + 1)}o$};
\draw (2*\c + 0.2, \c*0.6 - 0.3) node {$\w_{Nm_{2}}o$};
\draw (3*\c, 1.2*\c+0.18) node {$\w_{N(m_{2} + 0.5)}o$};
\draw (4*\c-0.14, \c*0.6 - 0.22) node {$\w_{N(m_{2} + 1)}o$};
\draw (5*\c + 0.2, \c*0.6 +0.2) node {$\w_{Nm_{3}}o$};
\draw (6*\c, -0.18) node {$\w_{N(m_{3} + 0.5)}o$};
\draw (7*\c-0.1, \c*0.6 + 0.18) node {$\w_{N(m_{3} + 1)}o$};
\draw (8*\c + 0.3, \c*0.6 - 0.2) node {$\w_{N \vartheta(M)}o$};

\end{tikzpicture}
\caption{A preliminary definition of pivotal loci}
\label{fig:scheme2}
\end{figure}

From these conditions, we can deduce small Gromov products among consecutive pivotal loci: $(w_{i-1}o, w_{i+1}o)_{w_{i}o}$ for each $i$. Note that this alone does not guarantee small Gromov products $(w_{i}o, w_{k}o)_{w_{j} o}$ for arbitrary $i<j<k$ that we need to sum up intermediate progresses: recall the theory of Mathieu-Sisto. This is remedied by the fact that small Gromov products are actually guided by long enough Schottky segments. Hence, each $d(w_{i-1} o, w_{i+1} o)$ is large enough and the Gromov inequality can deduce small arbitrary Gromov products. In $\T(\Sigma)$, we rely on the consequence of Rafi's theorems (cf. Section \ref{section:Teich}) and similarly deduce small arbitrary Gromov products.

The ideal situation is that all special steps can be hired as pivotal steps. If it is the case, the intermediate progresses are purely summed up and the overall progress becomes large enough. Unfortunately, there is always a small chance of the undesirable event: the probability that all special steps are pivotal times decays exponentially. Nonetheless, we want to `tolerate' such error and select large proportion of special steps that still performs the above task as pivotal times.

To further illustrate this idea, given pivotal times $\{m_{1}, \ldots, m_{k}\}$ for $(g_{1}, \ldots, g_{\vartheta(M)})$, let us determine the pivotal times for the path $(g_{1}, \ldots, g_{\vartheta(M+1)})$. One possible strategy is just adding $m_{k+1} = \vartheta(M)$ to the original set of pivotal times. Recall again the conditions for $m_{1}, \ldots, m_{k}$: \begin{itemize}
\item $(s_{i}^{-1} o, s'_{i}o)_{o} \le K$ for $i=1, \ldots, k$,
\item $(s'^{-1}_{i}o, w_{i}o)_{o} \le K$ for $i=1, \ldots, k$,
\item $(w_{i-1}^{-1} s'^{-1}_{i-1}o, s_{i} o)_{o} \le K$ $i=1, \ldots, k$. 
\end{itemize}
Fixing $w_{i}$ and $s_{i}'$, there are many other choices for each $s_{i}$ that satisfies the condition. In particular, the property of Schottky sets is designed so that at least $(\#S - 2)$ choices out of all choices are valid at each step. This process, fixing $w_{i}$ and $s_{i}'$ and modifying the choice of $s_{i}$ into another valid choice, is called \emph{pivoting}.

If, for example, the additional $s_{k+1}$, $s'_{k+1}$ and $w_{k+1}$ satisfy the above condition, then we can add it to the set of pivotal times. This already takes up large enough probability, at least $\left(\frac{\# S - 2}{\# S}\right)^{2}$. In case of failure, however, we do not wish to give up entire selection $\{m_{1}, \ldots, m_{k}\}$ but rather retain a portion that works for $(g_{1}, \ldots, g_{\vartheta(M+1)})$. For example, can we hope that the set $\{m_{1}, \ldots, m_{k-1}\}$ itself works for intermediate words $w_{0}$, $w_{1}$, $\ldots$, $w_{k-1} s_{k} s'_{k} w_{k} s_{k+1} s'_{k+1} w_{k+1}$? A priori, the final word depends on $s_{k+1} s'_{k+1}$ so this cannot be answered without altering $s_{k-1}'$: this is not what we want. We can however require the following condition: \begin{itemize}
\item $(s'^{-1}_{k-1}o, w_{k} o)_{o} \le K$,
\item $(w_{k-1}^{-1} s'^{-1}_{k-1} o, s_{k} o)_{o} \le K$,
\item $(s_{k}^{-1} o, s'_{k}w_{k}s_{k+1}s'_{k+1}w_{k+1} o)_{o} \le K$.
\end{itemize}
The first condition is already achieved by the fixed $s_{k-1}'$ from the assumption, and the latter two conditions can be achieved for any fixed $s_{k-1}'$, $w_{k-1}$, $s_{k+1}$, $s_{k+1}'$, $w_{k+1}$ by picking valid $s_{k}$ only. This has high chance so we have \[
\Prob\left[\{m_{1}, \ldots, m_{k-1}\} \,\textrm{works for}\, (g_{1}, \ldots, g_{\vartheta(M+1)}) | \{w_{i}, w'_{i}, s_{i}\}, \{s_{i}\}_{i \neq k}, s_{k}: \textrm{valid}\right] \ge\frac{\# S - 3}{\# S}.
\]
Note that the estimation is conditioned on each equivalence class of choices that are pivoted at $k$-th slot from each other. Summing them over bad choices of $s_{k+1}$, $s'_{k+1}$, we have \[
\Prob\left[\{m_{1}, \ldots, m_{k-1}\}, \{m_{1}, \ldots, m_{k}\}, \{m_{1}, \ldots, m_{k+1}\} \,\,\textrm{does not work}\right] \le \left[1-\left(\frac{\# S - 2}{\# S}\right)^{2}\right] \cdot \frac{3}{\# S}.
\]
Inductively, we deduce that the first $k-i$ slots (and possibly some more) among $\{m_{1}, \ldots, m_{k+1}\}$ can be employed as pivotal times for $(g_{1}, \ldots, g_{\vartheta(M+1)})$ except an exponentially decaying probability, whose decay rate depends on the size of $S$. In summary, we can guarantee almost definite increase of the number of pivots, as near as 1, by taking large enough Schottky set. For the precise definition that includes the modified conditions, see \cite{gouezel2021exp} or \cite{choi2021clt}.

Using small Gromov products among pivotal loci, one can show that $d(o, \w_{n} o)$ is bounded below by a multiple of $\#\{\textrm{pivots for}\,(g_{1}, \ldots, g_{n})\}$. Hence, we have established the definite progress of random walks outside an event of  exponentially decaying probability. In order to push this progress as close to the escape rate as we want, one should modify the decomposition \ref{eqn:schottkyMeas} and sandwich an auxiliary variable between Schottky steps. We refer the readers to \cite{gouezel2021exp} for details.

\section{Limit theorems II: translation lengths} \label{section:translation}

We now discuss the theory on translation lengths. In contrast with the case of displacements, where SLLN with the optimal moment condition was obtained at once, the first result on translation lengths was the following weak law of large numbers (WLLN).
\begin{thm}\label{thm:wllnTrans}
Let $X = \T(\Sigma)$ or $\C(\Sigma)$. Then there exists $L>0$ such that \[
\lim_{n} \Prob\left[\frac{1}{n} \tau_{X}(\w_{n}) \le L\right] =0.
\]
If $\mu$ further has finite first moment, then for any $\epsilon>0$ we have \[
\lim_{n} \Prob\left[\left|\frac{1}{n} \tau_{X}(\w_{n})- \lambda\right| > \epsilon\right] =0,
\]where $\lambda$ is the escape rate of the random walk.
\end{thm}
The WLLN for $\C(\Sigma)$ is proven by Maher-Tiozzo's theory, again due to the fact that the harmonic measure is atom-free. Its lifting to $\T(\Sigma)$ is due to Dahmani-Horbez's argument. Meanwhile, this convergence in probability alone is not enough to deduce the following SLLN.
\begin{thm}\label{thm:sllnTrans}
Let $X = \T(\Sigma)$ or $\C(\Sigma)$ and suppose that $\mu$ satisfies some moment condition. Then almost every random path $\w = (\w_{n})$ satisfies \[
\lim_{n} \frac{1}{n} \tau_{X}(\w_{n}) = \lambda.
\]
\end{thm}
In \cite{maher2018random}, Maher and Tiozzo discusses summable estimates of shadows along a direction for random walks on $\C(\Sigma)$ with bounded support. This is generalized to random walks with finite exponential moment, for which Boulanger, Mathieu, Sert and Sisto establish the same exponential decay of harmonic measure along the distance from the reference point \cite{boulanger2020large}. Moreover, as Dahmani and Horbez point out, Benoist-Quint's analog of Hsu-Robbins-Baum-Katz theorem gives summable estimates for random walks with finite second moment. This leads to Theorem \ref{thm:sllnTrans} under finite second moment, and Dahmani-Horbez's lifting copies this to $\T(\Sigma)$ as explained in Section \ref{section:Teich}. Here the lifting is possible whenever the random walk has finite first moment with respect to $d_{\T(\Sigma)}$, but it is the WLLN on $\C(\Sigma)$ that restricts the moment condition.

Before introducing Baik-Choi-Kim's theory in \cite{baik2021linear} for the SLLN under finite first moment condition, let us recall Maher-Tiozzo's strategy. In $\delta$-hyperbolic space such as $\mathcal{C}(\Sigma)$, the discrepancy $d(o, \w_{n}o) - \tau(\w_{n})$ is essentially correlated with the quantity $(\w_{n}^{-1} o, \w_{n} o)_{o}$. In particular, if $(\w_{n}^{-1} o, \w_{n} o)_{o}$ smaller than the half of $d(o, \w_{n} o)$ minus a constant, then the discrepancy is bounded by $(\w_{n}^{-1} o, \w_{n} o)_{o}$ plus a constant. 

In order to control $(\w_{n}^{-1} o, \w_{n} o)_{o}$, we now claim that the direction of $[o, \w_{n} o]$ ($[o, \w_{n}^{-1} o]$, resp.) is almost guided by $[o, g_{1} \cdots g_{\lfloor n/2\rfloor} o]$ ($[o, g_{n}^{-1} \cdots g_{\lfloor n/2 \rfloor + 1}^{-1} o]$, resp.). Given this claim, $(\w_{n}^{-1} o, \w_{n} o)_{o}$ is now correlated with the deviation $(g_{n}^{-1} \cdots g_{\lfloor n/2 \rfloor + 1}^{-1}o, g_{1} \cdots g_{\lfloor n/2\rfloor} o)$ between two \emph{independent} random paths. Actually, the claim itself also involves quantities of the same nature: $(\w_{n} o, \w_{\lfloor n/2 \rfloor} o)_{o}$ grows linearly almost surely if $d(o, \w_{\lfloor n/2} o)$ does grow linearly (which is true by the ergodic theorem) and $(o, \w_{n}o)_{\w_{\lfloor n/2\rfloor} o} = (g_{\lfloor n/2 \rfloor}^{-1} \cdots g_{1}^{-1} o, g_{\lfloor n/2 \rfloor+1} \cdots g_{n} o)_{o}$, the deviation between another pair of independent random paths, grows sublinearly.

Hence, it suffices to show that $\Prob[(\check{\w}_{n}o, \w_{n} o)_{o} \ge Kn]$ is summable for any $K>0$. For this Maher and Tiozzo conditions on each choice of $\check{\w}_{n}$ and regard $(\check{\w}_{n}o, \w_{n} o)_{o}$ as the deviation of a random segment $[o, \w_{n} o]$ from a fixed direction $[o, \check{\w}_{n} o]$. By Lemma 4.5 of \cite{benoist2016central} (together with the observation that $a$ can be chosen as any positive number for the item (4.8) in \cite{benoist2016central}), the probability is summable when $\mu$ has finite second moment and the conclusion follows. However, it is difficult to obtain summable estimates from weaker moment condition with this strategy, considering the Baum-Katz theorem.

The situation is more complicated in $\T(\Sigma)$ that lacks $\delta$-hyperbolicity. Here, summable estimates of $\Prob[(\check{\w}_{n}o, \w_{n} o)_{o} \ge Kn]$ is not enough and the deviation between paths should occur at special loci that enable arguments for $\delta$-hyperbolic spaces. Otherwise, one may give up controlling $\Prob[(\check{\w}_{n}o, \w_{n} o)_{o} \ge Kn]$ and pursue an argument not relying on the Borel-Cantelli lemma.

The theory of Baik-Choi-Kim falls into the latter case. Their strategy is to exploit the persistent joint of Maher and Tiozzo. Fixing suitable $L$, we say that a persistent joint arises at step $3kL$ if: \begin{itemize}
\item the steps $(g_{3(k-1)L+1}, \ldots, g_{3kL})$ constitute one of two Schottky-type paths,
\item the forward subpath $(\w_{3kL}o, \w_{3kL+1}o, \ldots)$ is contained in a shadow centered at $\w_{3kL}o$ viewed from $\w_{(3k-1)L} o$, and 
\item the backward subpath $(\ldots, \w_{3(k-1)L-1}o, \w_{3(k-1)L}o)$ is contained in a shadow centered at $\w_{(3k-2)L} o$. 
\end{itemize}
Both this construction and Gou{\"e}zel's construction aims to designate pivoting loci and derive almost sure phenomena. Nevertheless, the desired phenomena are different: Baik-Choi-Kim intend to guarantee large translation length by pivoting on the event of small translation length, while Gou{\"e}zel intends to guarantee definite progress from the prevalence of pivotal loci and one performs pivoting to establish this prevalence. 

Moreover, Baik-Choi-Kim's pivots entail technicalities that are not shared by Gou{\"e}zel's pivots. First, the prevalence of persistence joints is guaranteed outside events of summable probabilities. Another complication is that persistent joints are not independent variables. Nonetheless, persistent joints at different steps are linked by the ergodic shift and the subadditive ergodic theorem does guarantee the eventual prevalence of persistence joints for almost every path. We remark that Gou{\"e}zel's construction equally works with stronger implications.

We now explain how Baik-Choi-Kim achieved the almost sure linear growth of $\tau_{X}(\w_{n})$ without moment condition. We hope to declare an equivalence relation among random paths by pivoting: two paths are equivalent if they are identical except at the middle Schottky segment of the first $N$ persistent joints. Here comes one technical issue that the original and the pivoted path may not have the exactly same persistent joint steps. This is because persistent joints are random variables depending on entire $\w$ (not on finitely many steps near that joint) and single pivoting may alter entire distribution of persistent joints. To avoid this issue, Baik-Choi-Kim redefines pivots so that pivoting does not alter the pivot distribution. Moreover, according to their definition, persistent joints are incorporated in these pivots so the number of pivots also linearly grows almost surely.

Given this, the more pivots a path has, the smaller conditional probability that the path possess inside its equivalence class. We now observe that if a path $\w$ has small $\tau_{X}(\w_{n})$ and has enough number of early pivots within distance $\frac{1}{2}[d(o, \w_{n} o) - \tau_{X}(\w_{n})]$ from $o$, then the early pivoting $\w \mapsto \bar{\w}$ results in large $\tau_{X}(\bar{\w}_{n})$. Similar discussion holds for the late pivoting; in this case, enough number of late pivots within distance $\frac{1}{2}[d(o, \w_{n} o) - \tau_{X}(\w_{n})]$ from $\w_{n} o$ are needed. It remains to show that random paths either have enough number of early pivots near $o$, or have enough number of late pivots near $\w_{n}o$. This follows from distant allocation of pivotal loci: pivotal loci cannot be concentrated within distance $Ln$ for some suitable $L>0$, so those paths $\w$ with $\tau_{X}(\w_{n}) \le Ln$ necessarily fall into the above two categories.

The above argument can be improved if $\mu$ has finite first moment. For example, linearly growing number of pivots for $(g_{1}, \cdots, g_{n})$ should arise before $0.01n$ almost surely due to the subadditive ergodic theorem. Then the SLLN for displacements, another consequence of the ergodic theorem, asserts that linearly growing number of pivots appear within distance $0.02 \lambda n$ from $o$, where $\lambda$ is the escape rate. Thus, one can rely only on the early pivoting and bound the probability of $\{\tau_{X}(\w_{n}) \le d(o, \w_{n} o) - 0.04 \lambda n\}$. Since $d(o, \w_{n} o)/n$ also converges to $\lambda$, we obtain that $\limsup_{n} \tau_{X}(\w_{n})/n \ge 0.96 \lambda$ almost surely when $X = \C(\Sigma)$.

In the case of $\T(\Sigma)$, one should keep in mind that small $(\w_{n}^{-1} o, \w_{n}o)_{o}$ will not automatically imply small $(\w_{n}^{-m} o, \w_{n}^{k} o)_{o}$ for all $m, k>0$. Nonetheless, Baik-Choi-Kim exploits Rafi's results on thin triangles and fellow-traveling with thick ingredients and deduce the following fact. If the pivots are constructed with a $(K, K', \epsilon)$-Schottky set for sufficiently large $K'$, then each middle Schottky segment at the pivotal steps  are fellow traveling with some subsegment of $[o, \w_{n} o]$. Moreover, when a random path $\w$ satisfies $\frac{1}{2}[d(o, \w_{n} o) - \tau(\w_{n})] \ge d(o, z) + C$ for some pivotal locus $z$, then the directions of $[o, \check{\w}_{m} o]$ and $[o, \w_{n} o]$ near $z$ are guided by the same Schottky segment. If one pivots the path at $z$ by choosing different Schottky direction, then $[o, \check{\bar{\w}}_{m} o]$ and $[o, \bar{\w}_{n} o]$ deviate at $z$ and $[\check{\bar{\w}}_{m} o, \bar{\w}_{n} o]$ passes nearby $z$. This in turn implies that $\tau(\w_{n}) \ge d(o, \w_{n} o) - 2 d(o, z)$. Therefore, the SLLN in $\delta$-hyperbolic spaces is copied to $\T(\Sigma)$.

After Gou{\"e}zel and Baik-Choi-Kim's work, Choi tried to incorporate two notions of pivots in \cite{choi2021clt}. As a result, Choi explained how accurate the displacements and the translation lengths match from the prevalence of pivots except an exponentially decaying probability. This consequently implies the following deviation inequality.

\begin{prop}\label{prop:choiDeviate}
Suppose that $\mu$ has finite $p$-moment for some $p>0$ and let $q\le p$ be a nonnegative integer. Then there exists $K>0$ such that \[\begin{aligned}
\E \left[(\check{\w}_{m}o, \w_{m'}o)_{o}^{p+q}\right] <K+K e^{-m/K} (m' - m)^{q},\\ \quad \E \left[d(o, [\check{\w}_{m}o, \w_{m'}o])^{p+q}\right] < K + K e^{-m/K} (m' - m)^{q}
\end{aligned}
\]
for all $0 \le m \le m'$, respectively.
\end{prop}

In the special case $m = m'$, we obtain uniform control on  $\E[ (\check{\w}_{m}o, \w_{m'}o)_{o}^{2p}]$ from the finite $p$-moment of $\mu$. While Maher-Tiozzo's argument first fixes one of two random path and consider the deviation of the other path from that fixed direction, Choi performs pivoting on both random paths to make the estimate more effective and obtain exponent doubling.

In particular, when $\mu$ has finite $p$-moment for some $p>1/2$, the above estimate implies that \[
\Prob[(\check{\w}_{m} o, \w_{m} o)_{o} \ge Cm] \le \frac{\E[(\check{\w}_{m} o, \w_{m} o)_{o}^{2p}]}{(Cm)^{2p}} \le \frac{K}{(Cm)^{2p}}
\]
is summable for any $C>0$. Hence, Choi's result (together with Gou{\"e}zel's weak LDP from below) implies the SLLN for translation lengths in Gromov hyperbolic spaces including $\C(\Sigma)$ when $\mu$ has finite $p$-moment for some $p>1/2$. Nevertheless, this still requires a moderate moment condition; the SLLN for translation lengths without moment condition relies on the pivoting itself, as in Baik-Choi-Kim's argument. Choi deduces another consequence of the pivoting, the control of the discrepancy between displacements and translation lengths with greater precision: 

\begin{thm}
Suppose that $\mu$ has finite first moment. Then there exists a constant $K < \infty$ such that \begin{equation}
\limsup_{n} \frac{1}{\log n}| \tau(\w_{n}) - d(o, \w_{n} o)| < K
\end{equation}
for almost every $\w$.
\end{thm}

Meanwhile, the deviation inequality of Choi also turns out to be useful. One application is the improvement of the moment condition for geodesic tracking. Using the eventual version of Proposition \ref{prop:choiDeviate}, one can prove that sublinear geodesic tracking occurs in random walks with finite $(1/2)$-th moment. Moreover, a similar result for random walks with finite exponential moment implies logarithmic geodesic tracking. We note that logarithmic tracking was previously discussed on free groups, Gromov hyperbolic spaces and relatively hyperbolic spaces, the last two dealing with bounded support case (see  \cite{ledrappier2001free}, \cite{blachere2011harmonic}, \cite{sisto2017tracking}, \cite{maher2018random}).

Recall also that one can complete Mathieu-Sisto's approach to the CLT for displacements and translation lengths with this deviation inequality for $p=2$, hence achieving the optimal moment condition. Moreover, Choi established via explicit pivoting 
the converse of CLTs: the convergence of $\frac{1}{\sqrt{n}}[d(o, \w_{n} o) - c_{n}]$ or $\frac{1}{\sqrt{n}} [\tau(\w_{n}) - c_{n}]$ in law for some constant $c_{n}$ implies that the random walk has finite second moment. By adapting de Acosta's proof of the LIL for real-valued variables, Choi also establishes the LIL for displacements and translation lengths.

Let us now explain why the pivoting method is so effective. First, phenomena \emph{in probability} are correlated to certain probabilities that decay to zero, which ultimately relies on the non-atomness of the harmonic measure. This is due to the fact that $G$ is non-elementary: if a boundary point has the maximal atom, then all of its translations by $G$ should also have maximal atom and the boundary point should have finite orbit by $G$. This technique has been employed by many authors, including Woess \cite{woess1989boundaries}, Kaimanovich-Masur \cite{kaimanovich1996poisson}, Maher \cite{maher2011random} and Maher-Tiozzo \cite{maher2018random}.

Nonetheless, this is not enough deduce \emph{almost sure} phenomena and more accurate information is needed. Specifically, we need to elaborate the decaying rate of the harmonic measure corresponding to shadows, in terms of the distance of the shadows from the reference point. The first success was achieved by Maher, who deduced in \cite{maher2012exp} exponential decaying rate for random walks on $X = \C(\Sigma)$ with bounded support. We here explain a slight variation of Maher's argument.

Let us define the shadow $S_{x}(y, r)$ by the set $\{z : (y, z)_{x} \ge r\}$. Observe:
\begin{lem}\label{lem:shadow}
For $x, y, z \in G\cdot o$, $\xi \in X \cup \partial X$ and sufficiently large $R, R'>0$, we have the following: \begin{enumerate}
\item if $y \notin S_{x}(\xi, R')$ and $z \in S_{x}(\xi, R'+R)$, then $y \in S_{z}(o, R')$ and $z \in S_{y}(\xi, R')$.
\item $\nu(S_{x}(\xi, R)) := \Prob[ (\w_{n} x, \xi)_{x} \ge R\,\textrm{eventually}]\le 0.1$;
\item $H(S_{x}(\xi, R)) := \Prob[ (\w_{n} x, \xi)_{x} \ge R\,\textrm{at least once}]\le 0.12$;
\end{enumerate}
\end{lem}
The first item in fact holds for arbitrary $R$ and $R'$; it follows from the inequalities \[
(z, \xi)_{y}, (x, y)_{z} \ge (z, \xi)_{x} - (y, \xi)_{x},
\]
which are equivalent to the triangle inequality. The second item is due to the fact that $\nu$ is atom-free. For the last item, we should correlate the once-hitting event and the eventual event. Let $N_{0}$ be the first hitting time for $S_{x}(\xi, R)$, i.e., the earliest step at which $(\w_{n} x, \xi)_{o} \ge R$; this is a stopping time and the Markov property can be applied. Now with respect to any point $p \in S_{x}(\xi, R)$, we have $S_{x}(\xi, 0.5R)^{c} \subseteq S_{p}(x, 0.5R)$ by the second item. Since $\nu(S_{p}(o, 0.5R)) \le 0.1$ for sufficiently large $R$, we have $\nu [S_{x}(\xi, 0.5R) | \w_{N_{0}} x = p] \ge 0.9$ for each $p \in S_{x}(\xi, R)$. Consequently, we obtain $\nu[S_{x}(\xi, 0.5R)] \ge 0.9 H(S_{x}(\xi, R))$ and $H(S_{x}(\xi, R)) \le 0.12$.

Let us now fix $\xi \in X \cup \partial X$ and a sufficiently large number $R>0$, and estimate the hitting measure of $S_{o}(\xi, kR)$. We establish $k$ ``intermediate rivers" $R_{i}=S_{o}(\xi, \frac{3i-2}{3}R) \setminus S_{o}(\xi, \frac{3i-1}{3}R)$ that satisfy the following property:
\begin{enumerate}
\item each $R_{i}$ separates $X \setminus R_{i}$ into two part, $X_{i}^{+}$ and $X_{i}^{-}$, such that $d(X_{i}^{+}, X_{i}^{-}) > M$;
\item for each $i$, $R_{1}, \ldots, R_{i-1}$ are contained in $X_{i}^{-}$ and $R_{i+1}, \ldots, R_{k}$ are contained in $X_{i}^{+}$;
\item $R_{i+1}$ is contained in a shadow of distance $R/4$ with respect to any point in $R_{i}$.
\end{enumerate}

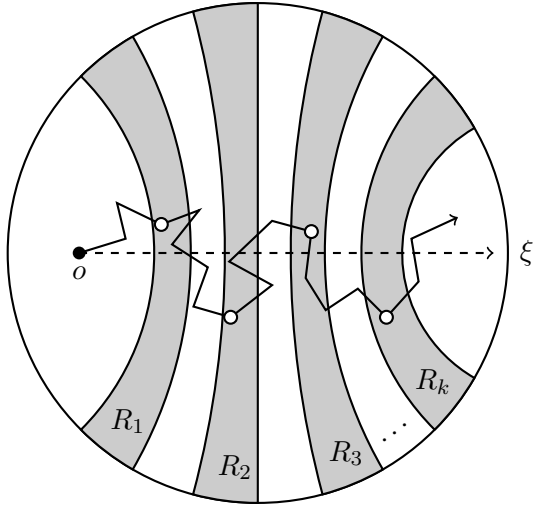
\begin{figure}[ht]
\begin{tikzpicture}[scale=0.95]
\draw[thick] (0, 0) circle (3.5);

\foreach \i in {2, 4}{
\draw[fill = black!20, thick] ({3.5 * cos(15*\i)}, {3.5* sin(15*\i)}) arc (90+15*\i : 270 - 15*\i : {3.5* tan(15*\i)}) arc (-15*\i : -15*\i - 15 : 3.5) arc (270 - 15* \i - 15 : 90+15*\i + 15 : {3.5* tan(15*\i + 15)}) arc (15*\i + 15 :15*\i : 3.5);
}

\draw[fill = black!20, thick] (0, 3.5) -- (0, -3.5) arc (-90 : -105 : 3.5) arc (-15:15:{3.5*tan(75)}) arc (105:90:3.5);

\draw[fill = black!20, thick] ({3.5*cos(120)}, {3.5*sin(120)}) arc (30 : -30: {3.5*tan(60)}) arc (-120 : -135 : 3.5) arc (-45:45:{3.5*tan(45)}) arc (135:120:3.5);

\draw[thick, ->] (-2.5, 0) -- (-1.85, 0.2) -- (-1.97, 0.7) -- (-1.35, 0.4) -- (-0.82, 0.6) -- (-1.2,0.12) -- (-0.7, -0.2) -- (-0.9, -0.75) --  (-0.38, -0.9) -- (0.2, -0.45) -- (-0.4, -0.12) -- (0.2, 0.45) -- (0.75, 0.3) -- (0.67, -0.35) -- (0.95, -0.8) -- (1.4, -0.5) -- (1.8, -0.9) -- (2.25, -0.4) -- (2.15, 0.2) -- (2.8, 0.5);

\draw[thick, fill=black!0] (-1.35, 0.4) circle (0.09);
\draw[thick, fill=black!0] (-0.38, -0.9) circle (0.09);
\draw[thick, fill=black!0]  (0.75, 0.3) circle (0.09);
\draw[thick, fill=black!0]  (1.8, -0.9) circle (0.09);

\draw[thick, dashed, ->] (-2.5, 0) -- (3.3, 0);

\draw (-2.5, -0.28) node {$o$};
\fill (-2.5, 0) circle (0.09);
\draw (3.765, 0) node {$\xi$};
\draw (-1.82, -2.35) node {$R_{1}$};
\draw (-0.32, -3) node {$R_{2}$};
\draw (1.235, -2.8) node {$R_{3}$};
\draw[shift={(1.95, -2.5)}] (0, 0) node[rotate=37.5] {$\cdots$};
\draw (2.45, -1.8) node {$R_{k}$};

\end{tikzpicture}
\caption{Schematics of rivers $R_{i}$ in Maher's argument. Each white dots represents the sample locus at $N_{i}(\w)$.}
\end{figure}

Due to properties (1), (2) and the fact that each random path consists of bounded steps, it is necessary to enter each $R_{i}$ at least once to reach beyond $R_{k}$. This motivates us to consider the first hitting time $N_{i}(\w)$ at which step $\w_{n} o$ first enters $R_{i}$ and define $E_{i} := \{\w : N_{1}(\w), \ldots, N_{i}(\w) < \infty\}$. In order to calculate $\Prob[E_{i+1}|E_{i}]$, we condition on each choice of $\w$ until $N_{i}$ and fix $p = \w_{N_{i}}o \in R_{i}$. As $N_{i}(\w)$ is a stopping time, one can then apply the Markov property for estimation. Namely, $R_{i+1} \subseteq S_{o}(\xi, \frac{3i+1}{3}R) \subseteq S_{p}(\xi, 2R/3)$ has hitting measure at most 0.12, and we have \[\begin{aligned}
\Prob[E_{i+1}] &= \sum_{a_{1}, \ldots, a_{n}\in G} H(R_{i+1}) \cdot \Prob[\w: N_{i}(\w) = n, g_{i} = a_{i}\,\textrm{for each}\, i=1, \ldots, n]\\
&\le \sum_{a_{1}, \ldots, a_{n}\in G} H(S_{a_{1} \cdots a_{n}o}(\xi, 2R/3)) \cdot \Prob[\w: N_{i}(\w) = n, g_{i} = a_{i}\,\textrm{for each}\, i=1, \ldots, n]\\
&\le 0.12 \sum_{a_{1}, \ldots, a_{n}\in G}\Prob[\w: N_{i}(\w) = n, g_{i} = a_{i}\,\textrm{for each}\, i=1, \ldots, n] = 0.12 \Prob[E_{i}].
\end{aligned}
\]
This implies $H(S_{o}(\xi, kR)) \le \Prob[E_{k}] \le 0.12^{k}$ as desired.

The role of Lemma \ref{lem:shadow} is to correlate the probability of progress in a specific direction with the probability of the rest. Note that this process does not require moment condition. However, in order to correlate those probabilities with the distance, one need to quotient out the path space into measurable equivalence classes at regular distances. This is realized as hitting times, which crucially depends on the boundedness of each step so that no path jumps over and skips any river. We remark that random walks with finite exponential moment exhibit similar behavior. Although the hitting time for each river is not exactly realized, the probability of the error case that one jumps over $n$ river decays exponentially and we have similar exponential decay of the harmonic measure (cf. Corollary 2.13, \cite{boulanger2020large}). 

Meanwhile, Benoist-Quint's martingale version of Hsu-Robbins-Baum-Katz theorem can deal with measures with finite $p$-moment, beyond those with bounded support. This method aims to estimate the concentration of the cocycles $\sigma(g, x)$ (measured with $\mu^{\ast n}$) near the average $\lambda n$. In this perspective, one adds up $n$ martingale differences $\sigma_{0}(g_{n}^{-1}, g_{n-1}^{-1} \cdots g_{1}^{-1} o)$, which is a balanced version of $d(x, \w_{n} o) - d(x, \w_{n-1} o)$ and is bounded by $2d(o, g_{n} o)$. Since this step is $L^{p}$-bounded, one obtains an $L^{p-2}$-convergence rate: for each $\epsilon>0$, there exist constants $D_{n}$ such that $\sum_{n} n^{p-2} D_{n} < \infty$ and \[
\Prob\left[\w : \lambda - \epsilon \le \sigma(\w_{n}^{-1}, x) \le \lambda+\epsilon \right]\le D_{n}.
\]

The advantage of this method is that it applies to $L^{p}$-integrable cocycles on arbitrary compact metric space, where the $p$-th moments of the steps are uniformly bounded but not summable: it becomes summable only after modulating by order 2. In our setting of Gromov hyperbolic spaces of Teichm{\"u}ller spaces, however, one can expect further efficiency from below since the $p$-th moment of the steps are not only bounded but exponentially decaying. More precisely, one does not observe $d(x, \w_{n} o) - d(x, o)$ but observes its counterpart $(\w_{n} o, x)_{o}$: this conversion requires Gromov hyperbolicity or its analogy on Teichm{\"u}ller spaces. 

This alternative strategy is also pursued in \cite{benoist2016central}, beginning from the spectral gap of amenable groups acting on compact spaces and the exponential growth. Recall that another approach to this part of the argument, suggest in \cite{choi2021clt}, removes the cocompactness assumption. Given this preliminary estimates, the final step is to bound $\lim_{n} (\w_{n} o, x)_{o}^{p}$ with the sum of $d(o, g_{k+1}o)^{p} 1_{\lim (x, \w_{n} o) \ge d(o, \w_{k} o)}$ (when $0<p<1$) or $2^{p} [d(o, g_{k+1}o)^{p} + d(o, \w_{k}o)^{p-1} d(o, g_{k+1}o)]1_{\lim (x, \w_{n} o) \ge d(o, \w_{k} o)}$ (when $p\ge 1)$ and control each expectation. In a plain language, this counts the contribution of each step of the form $d(o, g_{k+1}o)^{p}$ or $d(o, \w_{k}o)^{p-1} d(o, g_{k+1}o)$ only when the progress of the random path is toward $x$, whose probability decays exponentially and results in summable contribution to the $p$-th moment.

The description so far of Maher's and Benoist-Quint's theories are by no means complete; for fuller analysis, including the martingale version of Hsu-Robbins-Baum-Katz theorem, see \cite{maher2012exp}, \cite{benoist2016central} and \cite{benoist2016linear}.

To sum up, the above strategies estimate the decay rate of the harmonic measure via packing random paths into effective and ineffective cases, by referring to the intermediate steps, and count the effective cases only. This philosophy is maximized in the notion of pivots in Gou{\"e}zel's, Baik-Choi-Kim's and Choi's work. Each equivalent class of the same pivots consists mostly of the desirable paths and a small portion of undesirable paths; their probability can be compared by pivoting and can be summed using the Markov property. The notion of pivots also fit into the realm of random walks with infinite support; moment conditions are not necessary for the punctual appearance of pivots, and are used only to synchronize the time and distance progress of pivots.

\section{Distance and counting}

So far, we have discussed various methods to study random mapping classes that arise from random walks, especially with respect to their action on $\C(\Sigma)$ or $\T(\Sigma)$. This philosophy essentially differs from studying the random mapping classes on $\Mod(\Sigma)$ itself, since $\Mod(\Sigma)$ equipped with the word metric is not quasi-isomorphic to $\C(\Sigma)$ nor $\T(\Sigma)$. Hence, counting elements in $\Mod(\Sigma)$ with respect to the word metric becomes a separate problem.

One possible solution is to use properties of the action of $\Mod(\Sigma)$ on $\T(\Sigma)$ or $\C(\Sigma)$ beyond non-elementariness. In this direction, we mentioned that Mathieu and Sisto exploited the acylindrical hyperbolicity of $\Mod(\Sigma)$ to bring the aforementioned results (including definite progress, CLT, etc.) on $\C(\Sigma)$ to $\Mod(\Sigma)$ or $\T(\Sigma)$ \cite{mathieu2020deviation}.

Another solution is to realize a Markov process on the group itself. Here is used the automatic structure of groups, first hinted by Cannon \cite{cannon1984combinatorial} and later formulated by W. Thurston. For general reference, see \cite{epstein1992word}. An automatic structure of a group models (quasi-)geodesics on the group with paths on a directed graph. By considering a Markov process on this graph, we can utilize the techniques for random walks to describe the asymptotic behavior in the counting setting. In particular, if the graph possesses suitable hyperbolicity (e.g. exponential growth, independent directions, etc.), then the counting problem (guided by the Patterson-Sullivan measure) mingle with the random walk theory (guided by the harmonic measure). One can also interpret pivots as a partial realization of automatic structure, by recording pivotal times (as if they represent specific vertex on the graph structure) and pivoting the choices at pivots (as if we distinguish cone types).

Notable examples of (geodesic) automatic groups include hyperbolic groups, relatively hyperbolic groups, right-angled Artin/Coxeter groups and many more. In particular, hyperbolic groups have geodesic automatic structure with respect to any finite generating set, allowing the WLLN \cite{gekhtman2018counting} and CLT \cite{gekhtman2020clt} for for displacements and translation lengths. Still, the theory is not applicable for the entire mapping class group at the moment: although mapping class groups have quasi-geodesic automatic structure, it is not known whether they have geodesic automatic structure. Nonetheless, if the generating set is nicely populated with some Schottky set, one can partially realize this strategy on weakly hyperbolic groups and $\Mod(\Sigma)$ that lack geodesic automatic structure. This will be explained further in the forthcoming preprint \cite{choi2021pA} of the second author.

\section{Further directions}

We have discussed the random walks on $\Mod(\Sigma)$ in different perspectives. Several questions arise from the difference among groups and spaces. First, it is known that random walks on hyperbolic groups also satisfy local limit 
theorem \cite{gouezel2014local}. The ingredient of Gou{\"e}zel's argument that depends on the Gromov hyperbolicity is to establish Ancona's inequality. Considering the parallel theory of pivoting on Gromov hyperbolic spaces and $\T(\Sigma)$, one might hope a similar result on $\T(\Sigma)$.

There is a twin notion for $\Mod(\Sigma)$ acting on $\T(\Sigma)$ and $\C(\Sigma)$, namely, the outer automorphism group $Out(F_{n})$ that acts on the Culler-Vogtmann Outer space and the complex of free factors. As in the case of $\Mod(\Sigma)$, the dynamical property of $\varphi \in Out(F_{n})$ is revealed by its action on the Outer space, and the complex of free factors is often chosen as a detour since it is Gromov hyperbolic. We expect that our theory for random walks on $\Mod(\Sigma)$ is almost exactly transcribed into the one on $Out(F_{n})$. In fact, some progress were already made by the work of Horbez \cite{horbez2018central} and Dahmani-Horbez \cite{dahmani2018spectral}.

Finally, despite partial achievements, the complete Patterson-Sullivan theory on $\Mod(\Sigma)$ is not attained yet. Once achieved, this will serve as another perspective for the counting problem in $\Mod(\Sigma)$. For instance, see Gekhtman's analysis on the stable type of mapping class groups \cite{gekhtman2013stable}.

\bibliographystyle{alpha} 
\bibliography{sullivan}

\end{document}